\documentclass[11pt,a4paper]{article}

\usepackage[T1]{fontenc}
\usepackage{lmodern}
\usepackage{microtype}

\usepackage{amsmath,amssymb,mathtools,amsthm}
\usepackage[all,cmtip]{xy}
\usepackage{graphicx}
\usepackage{float}
\usepackage{xcolor}
\usepackage{tikz}
\usepackage{tikz-cd}
\usepackage[a4paper,margin=1in]{geometry}

\usepackage[backend=biber,style=numeric]{biblatex}

\usepackage[unicode=true,pdfusetitle,bookmarks=true,bookmarksnumbered=false,
  bookmarksopen=false,breaklinks=true,pdfborder={0 0 0},backref=false,
  colorlinks=true]{hyperref}
\usepackage{cleveref}
\usepackage{todonotes}

\newtheorem{theorem}{Theorem}[section]
\newtheorem{lemma}[theorem]{Lemma}
\newtheorem{proposition}[theorem]{Proposition}
\newtheorem{corollary}[theorem]{Corollary}

\theoremstyle{remark}
\newtheorem{remark}[theorem]{Remark}

\AddToHook{env/theorem/begin}{\crefalias{theorem}{theorem}}
\AddToHook{env/lemma/begin}{\crefalias{theorem}{lemma}}
\AddToHook{env/proposition/begin}{\crefalias{theorem}{proposition}}
\AddToHook{env/corollary/begin}{\crefalias{theorem}{corollary}}
\AddToHook{env/remark/begin}{\crefalias{theorem}{remark}}

\makeatletter
\AddToHook{cmd/thmhead/before}{\MakeLinkTarget*{\@currentHref.heading}}
\makeatother

\colorlet{linkColor}{blue!55!black}
\colorlet{citationColor}{red!55!black}
\colorlet{sourceFill}{blue}
\colorlet{sourceBoundary}{blue!70!black}
\colorlet{targetFill}{orange!15}
\colorlet{targetBoundary}{orange!80!black}
\colorlet{arithmeticArrowColor}{blue!55!black}
\colorlet{figureAxisColor}{black!80}
\colorlet{figurePointColor}{black!88}
\colorlet{figureGuideColor}{black!65}
\colorlet{figureHollowFill}{white}
\definecolor{torsionLevelSeven}{HTML}{EFF3FF}
\definecolor{torsionLevelEight}{HTML}{BDD7E7}

\hypersetup{linkcolor=linkColor,urlcolor=linkColor,citecolor=citationColor}

\newcommand{\bbR}{\mathbb{R}}
\newcommand{\bbC}{\mathbb{C}}
\newcommand{\bbZ}{\mathbb{Z}}
\newcommand{\bbN}{\mathbb{N}}
\newcommand{\bbQ}{\mathbb{Q}}
\newcommand{\bbF}{\mathbb{F}}

\newcommand{\Cn}{\mathbb{C}^{n}}
\newcommand{\cA}{\mathcal{A}}

\newcommand{\cL}{\mathcal{L}}

\newcommand{\SH}{\operatorname{SH}}

\newcommand{\im}{\operatorname{im}}

\title{Torsion phenomena in symplectic capacities via Clarke duality}
\author{Arye Deutsch\thanks{Einstein Institute of Mathematics, Hebrew University of Jerusalem, Israel. \href{mailto:arye.deutsch@mail.huji.ac.il}{arye.deutsch@mail.huji.ac.il}}}
\date{}

\begin{document}

\maketitle

\begin{abstract}
We study Gutt--Hutchings capacities with finite cyclic coefficients
and show that torsion in positive $S^1$-equivariant symplectic homology
detects new obstructions to symplectic embeddings between convex toric
domains.  Prime-field coefficients give no additional obstructions
in this class, while the torsion responsible for the new obstructions
already appears over $\mathbb Z/p^{p+2}$ for every prime $p$.
We prove a stabilization theorem for these capacities on convex
bodies with arbitrary commutative coefficients, yielding a sharp
stabilized embedding obstruction in an explicit example.
Clarke duality makes these phenomena accessible to computation.
We extend the filtered equivariant Clarke--Floer comparison to
arbitrary commutative coefficient rings and compute the capacities
using models of dual action sublevels built from representation
spheres.  The circle-action weights encode orbit covering
multiplicities, whose products govern the new obstructions.
\end{abstract}

\tableofcontents
\section{Introduction}

Symplectic capacities distinguish symplectic embeddings from
volume-preserving embeddings, as illustrated by Gromov's
non-squeezing theorem.  They give obstructions beyond volume:
a symplectic embedding $X\hookrightarrow Y$ requires the capacity
of $X$ to be at most that of $Y$.

The capacities considered here come from positive $S^1$-equivariant
symplectic homology.  They record the action levels at which
distinguished homology classes appear.  The usual Gutt--Hutchings
capacities use rational coefficients \cite{gutt2018symplectic}.
In this paper we let $R$ vary over commutative unital rings, write
$c_k^R$ for the resulting capacities, and ask whether this change reveals
embedding obstructions invisible to the complete rational sequence.
The new obstructions arise not from introducing new moduli spaces,
but from changing the coefficient ring in which the parametrized
Floer cylinders of the Gutt--Hutchings construction are counted.
For star-shaped domains $X,Y\subset\bbC^n$, these capacities satisfy
symplectic monotonicity, as in
\cite[Theorem~1.1 and Remark~1.3]{gutt2018symplectic}:
\[
 X\hookrightarrow Y\quad\Longrightarrow\quad
 c_k^R(X)\leq c_k^R(Y)\qquad(k\geq1).
\]

The dependence on coefficients has also been studied over fields.
\c{C}ineli, Ginzburg, and G\"urel
\cite[Theorem~A]{cineli2026filtered} show that the sequence over a
field of positive characteristic is bounded above, in contrast to
the rational sequence.  For irrational ellipsoids they compute it
as the rational sequence truncated at the action of the first
$p$-fold iterate \cite[Example~1.2]{cineli2026filtered}.
They ask whether one can compute these invariants on convex or
concave toric domains, and how much of the ellipsoid behavior
persists more generally \cite[Section~1.2.1]{cineli2026filtered}.

Changing the capacity values does not by itself give additional
embedding obstructions.  For convex toric domains, we compute the
prime-field capacities and find that all their inequalities are
already implied by the rational ones.

\begin{proposition}
\label{prop:prime-field-toric-capacities}
For every convex toric domain $X$, every prime $p$, and every $k\geq1$,
\begin{equation}
 \label{eq:prime-field-toric-capacities}
 c_k^{\bbF_p}(X)
       =\min\{c_k^{\bbQ}(X),p\,c_1^{\bbQ}(X)\}.
\end{equation}
Consequently, for convex toric domains $X$ and $Y$ of the same
dimension, the inequalities
$c_j^{\bbQ}(X)\leq c_j^{\bbQ}(Y)$ for every $j\geq1$ imply
$c_k^{\bbF_p}(X)\leq c_k^{\bbF_p}(Y)$ for every prime $p$ and
every $k\geq1$.
\end{proposition}

We prove this in \cref{sec:cyclic-coefficient-capacities}, combining
the supporting-ellipsoid bound with the low-action coefficient
comparison of \c{C}ineli, Ginzburg, and G\"urel.
The formula answers their calculation question for convex toric
domains: the sequence is eventually constant at $p\,c_1^{\bbQ}(X)$,
which for smooth domains is the least action of a $p$-fold iterate.

Prime-field coefficients therefore give no additional embedding
obstructions between convex toric domains.  For embeddings between
ellipsoids, even finite cyclic coefficients give no additional
obstructions; see \cref{thm:no-new-cyclic-ellipsoid-obstruction}.
For more general convex toric domains, however, finite cyclic
coefficients can reveal obstructions invisible to the complete
rational sequence.  We exhibit such a pair using coefficients in
$\bbZ/243$.

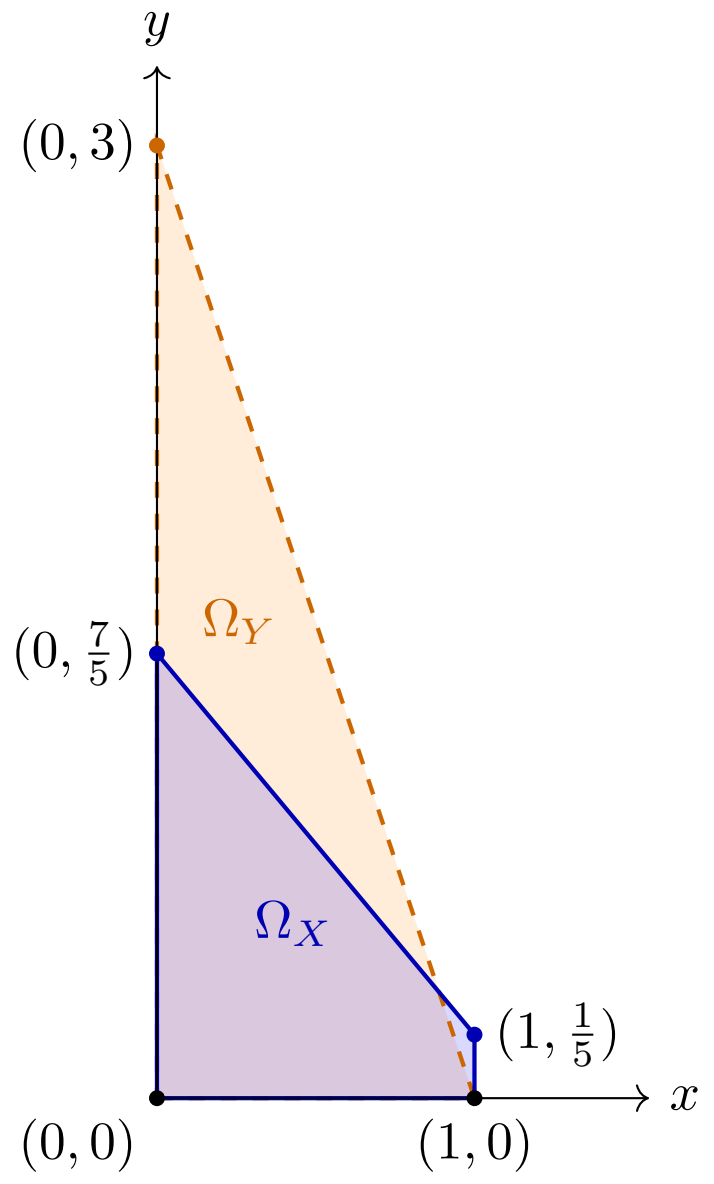
\begin{figure}[htbp]
    \centering
    \begin{tikzpicture}[x=2.1cm,y=2.1cm,
        every node/.style={font=\small}]
        \fill[targetFill] (0,0) -- (1,0) -- (0,3) -- cycle;
        \fill[sourceFill,fill opacity=0.15]
            (0,0) -- (1,0) -- (1,{1/5}) -- (0,{7/5}) -- cycle;

        \draw[targetBoundary,thick,dashed]
            (0,0) -- (1,0) -- (0,3) -- cycle;
        \draw[sourceBoundary,thick]
            (0,0) -- (1,0) -- (1,{1/5}) -- (0,{7/5}) -- cycle;

        \draw[->] (0,0) -- (1.55,0) node[right] {$x$};
        \draw[->] (0,0) -- (0,3.25) node[above] {$y$};

        \fill[targetBoundary] (0,3) circle (1.5pt);
        \node[left] at (0,3) {$(0,3)$};
        \fill[sourceBoundary] (0,{7/5}) circle (1.5pt);
        \node[left] at (0,{7/5}) {$(0,\frac75)$};
        \fill[sourceBoundary] (1,{1/5}) circle (1.5pt);
        \node[right] at (1,{1/5}) {$(1,\frac15)$};
        \fill (1,0) circle (1.5pt);
        \node[below] at (1,0) {$(1,0)$};
        \fill (0,0) circle (1.5pt);
        \node[below left] at (0,0) {$(0,0)$};

        \node[sourceBoundary] at (0.43,0.55) {$\Omega_X$};
        \node[targetBoundary] at (0.26,1.5) {$\Omega_Y$};
    \end{tikzpicture}
    \caption{Overlaid moment regions of $X$ (blue, solid boundary) and
    $Y$ (orange, dashed boundary).}
    \label{fig:243-moment-regions}
\end{figure}

Write
\[
 \mu(z_1,\ldots,z_n)=\pi(|z_1|^2,\ldots,|z_n|^2),
 \qquad X_\Omega=\mu^{-1}(\Omega),
\]
and let
\[
 E(a_1,\ldots,a_n)=
 \left\{z\in\bbC^n\ \middle|\
           \sum_i\frac{\pi|z_i|^2}{a_i}\leq1\right\}.
\]
Our source is $X=X_{\Omega_X}$ and our target is $Y=E(1,3)$, where
\[
 \Omega_X=
 \left\{(x,y)\in\bbR_{\geq0}^2\ \middle|\
             x\leq1,\quad 6x+5y\leq7\right\}.
\]
The moment regions are shown in \cref{fig:243-moment-regions}.

% Moment regions from TorsionObstruction243.tex.
% A single coordinate system with equal x and y units.

\begin{theorem}
\label{thm:243-obstruction}
The capacity with coefficients in $\bbZ/243$ satisfies
\[
 c_{15}^{\bbZ/243}(X)=\frac{48}{5}>9=c_{15}^{\bbZ/243}(Y),
\]
whereas the rational Gutt--Hutchings and ECH capacities satisfy
\[
 c_j^{\bbQ}(X)\leq c_j^{\bbQ}(Y),\qquad
 c_j^{\mathrm{ECH}}(X)\leq c_j^{\mathrm{ECH}}(Y)
 \qquad\text{for every }j\geq1.
\]
Consequently, there is no symplectic embedding
\[
 \operatorname{int}(X)\hookrightarrow\operatorname{int}(Y).
\]
\end{theorem}

The capacity values are computed in \cref{sec:243-domains}.

The same obstruction, still invisible to the rational and ECH
capacities, holds for sufficiently close smooth strongly convex inner
approximations of $X$, by continuity and monotonicity of capacities.

Volume also gives no obstruction:
$\operatorname{vol}(X)=4/5<3/2=\operatorname{vol}(Y)$.
In dimension four, Hutchings' refined ECH criterion
\cite[Theorem~1.19]{hutchings2016beyond} also yields the sharp
bound $\lambda\geq16/15$ for embeddings into $E(\lambda,3\lambda)$,
corresponding to $N=0$ in \cref{cor:stabilized-243} below.
To obtain an obstruction after stabilization, we prove that
sufficiently large disk factors preserve these capacities for every
convex body and every coefficient ring.

Write
\[
 D(A)=\{w\in\bbC:\pi|w|^2\leq A\}
\]
for the disk of area $A$.

\begin{theorem}[Stabilization]
\label{thm:capacity-disk-stabilization}
Let $C\subset\bbC^n$ be a compact convex body containing the origin
in its interior.  For a commutative unital ring $R$, $k\geq1$, and
$A_1,\ldots,A_N>c_k^R(C)$, one has
\[
 c_k^R\bigl(C\times D(A_1)\times\cdots\times D(A_N)\bigr)=c_k^R(C).
\]
\end{theorem}

Applying \cref{thm:capacity-disk-stabilization} to the capacity values
in \cref{thm:243-obstruction} gives the following sharp embedding threshold.

\begin{corollary}
\label{cor:stabilized-243}
For every integer $N\geq0$ and every $\lambda>0$, there is a symplectic
embedding
\[
 \operatorname{int}(X)\times\bbC^N
       \hookrightarrow
 \operatorname{int}(E(\lambda,3\lambda))\times\bbC^N
\]
if and only if $\lambda\geq16/15$.
\end{corollary}

\begin{figure}[htbp]
    \centering
    \resizebox{\linewidth}{!}{%
    \begin{tikzpicture}[x=2.20cm,y=1.65cm,>=stealth]
        \draw[->] (2.4,0.45) -- (11.55,0.45);
        \draw[->] (2.4,0.45) -- (2.4,5.45);
        \node at (7.0,-0.05) {homological degree};
        \node[rotate=90] at (1.75,3.0) {action};
        \foreach \x/\d in {3/2,4/3,5/4,6/5,7/6,8/7,9/8,10/9,11/10} {
            \draw (\x,0.39) -- (\x,0.51);
            \node[below] at (\x,0.39) {$\d$};
        }
        \foreach \y/\lab in {1/A,2/{2A},3/{3A},4/{4A},5/{5A}} {
            \draw (2.34,\y) -- (2.46,\y);
            \node[left] at (2.34,\y) {$\lab$};
        }
        \node (r1u0q) at (3,1) {$x_1$};
        \node (r1u0p) at (4,1) {$y_1$};
        \node (r1u1q) at (5,1) {$ux_1$};
        \node (r1u1p) at (6,1) {$uy_1$};
        \node (r1u2q) at (7,1) {$u^{2}x_1$};
        \node (r1u2p) at (8,1) {$u^{2}y_1$};
        \node (r1u3q) at (9,1) {$u^{3}x_1$};
        \node (r1u3p) at (10,1) {$u^{3}y_1$};
        \node[text=sourceBoundary] (r1u4q) at (11,1) {$u^{4}x_1$};
        \node (r2u0q) at (5,2) {$x_2$};
        \node (r2u0p) at (6,2) {$y_2$};
        \node (r2u1q) at (7,2) {$ux_2$};
        \node (r2u1p) at (8,2) {$uy_2$};
        \node (r2u2q) at (9,2) {$u^{2}x_2$};
        \node (r2u2p) at (10,2) {$u^{2}y_2$};
        \node[text=sourceBoundary] (r2u3q) at (11,2) {$u^{3}x_2$};
        \node (r3u0q) at (7,3) {$x_3$};
        \node (r3u0p) at (8,3) {$y_3$};
        \node (r3u1q) at (9,3) {$ux_3$};
        \node (r3u1p) at (10,3) {$uy_3$};
        \node[text=sourceBoundary] (r3u2q) at (11,3) {$u^{2}x_3$};
        \node (r4u0q) at (9,4) {$x_4$};
        \node (r4u0p) at (10,4) {$y_4$};
        \node[text=sourceBoundary] (r4u1q) at (11,4) {$ux_4$};
        \node[text=targetBoundary] (r5u0q) at (11,5) {$x_5$};
        \draw[->] (r1u1q) -- node[above] {$1$} (r1u0p);
        \draw[->] (r1u2q) -- node[above] {$1$} (r1u1p);
        \draw[->] (r1u3q) -- node[above] {$1$} (r1u2p);
        \draw[->] (r1u4q) -- node[above] {$1$} (r1u3p);
        \draw[->] (r2u0q) -- node[midway,above,inner sep=1.5pt,text=arithmeticArrowColor] {$1$} (r1u0p);
        \draw[->] (r2u1q) -- node[above] {$2$} (r2u0p);
        \draw[->] (r2u1q) -- node[midway,above,inner sep=1.5pt,text=arithmeticArrowColor] {$1$} (r1u1p);
        \draw[->] (r2u2q) -- node[above] {$2$} (r2u1p);
        \draw[->] (r2u2q) -- node[midway,above,inner sep=1.5pt,text=arithmeticArrowColor] {$1$} (r1u2p);
        \draw[->] (r2u3q) -- node[above] {$2$} (r2u2p);
        \draw[->] (r2u3q) -- node[midway,above,inner sep=1.5pt,text=arithmeticArrowColor] {$1$} (r1u3p);
        \draw[->] (r3u0q) -- node[midway,above,inner sep=1.5pt,text=arithmeticArrowColor] {$1$} (r2u0p);
        \draw[->] (r3u1q) -- node[above] {$3$} (r3u0p);
        \draw[->] (r3u1q) -- node[midway,above,inner sep=1.5pt,text=arithmeticArrowColor] {$1$} (r2u1p);
        \draw[->] (r3u2q) -- node[above] {$3$} (r3u1p);
        \draw[->] (r3u2q) -- node[midway,above,inner sep=1.5pt,text=arithmeticArrowColor] {$1$} (r2u2p);
        \draw[->] (r4u0q) -- node[midway,above,inner sep=1.5pt,text=arithmeticArrowColor] {$1$} (r3u0p);
        \draw[->] (r4u1q) -- node[above] {$4$} (r4u0p);
        \draw[->] (r4u1q) -- node[midway,above,inner sep=1.5pt,text=arithmeticArrowColor] {$1$} (r3u1p);
        \draw[->] (r5u0q) -- node[midway,above,inner sep=1.5pt,text=arithmeticArrowColor] {$1$} (r4u0p);
        % Coefficients are aligned with the corresponding terms of Z_5.
        \node at (11,5.95) {$Z_5$};
        \draw[->] (11,5.72) -- (11,5.25);
        \node[font=\small] at (11.85,5.60) {coefficients};
        \draw (11.62,5.22) -- (11.52,5.22) -- (11.52,0.78) -- (11.62,0.78);
        \draw (12.08,5.22) -- (12.18,5.22) -- (12.18,0.78) -- (12.08,0.78);
        \foreach \y/\c in {1/1,2/{-1},3/2,4/{-6}} {
            \node[text=sourceBoundary] at (11.85,\y) {$\c$};
        }
        \node[text=targetBoundary] at (11.85,5) {$24$};
    \end{tikzpicture}%
    }
    \caption{The positive $S^1$-equivariant complex of the disk $D(A)$.
    The differential is given by Zhao's disk computation
    \cite[Section~8.1]{zhao2019periodic}, expressed here in the positive
    equivariant homological complex.
    The coefficients beside the rightmost column specify the cycle
    $Z_5$ defining the fifth capacity; its blue terms lie at action
    at most $4A$, and its orange term lies at action $5A$.
    The generators $x_j,y_j$ correspond to the $j$-fold cover,
    with $|x_j|=2j$, $|y_j|=2j+1$, and $|u|=2$.}
    \label{fig:disk-arithmetic-complex}
\end{figure}
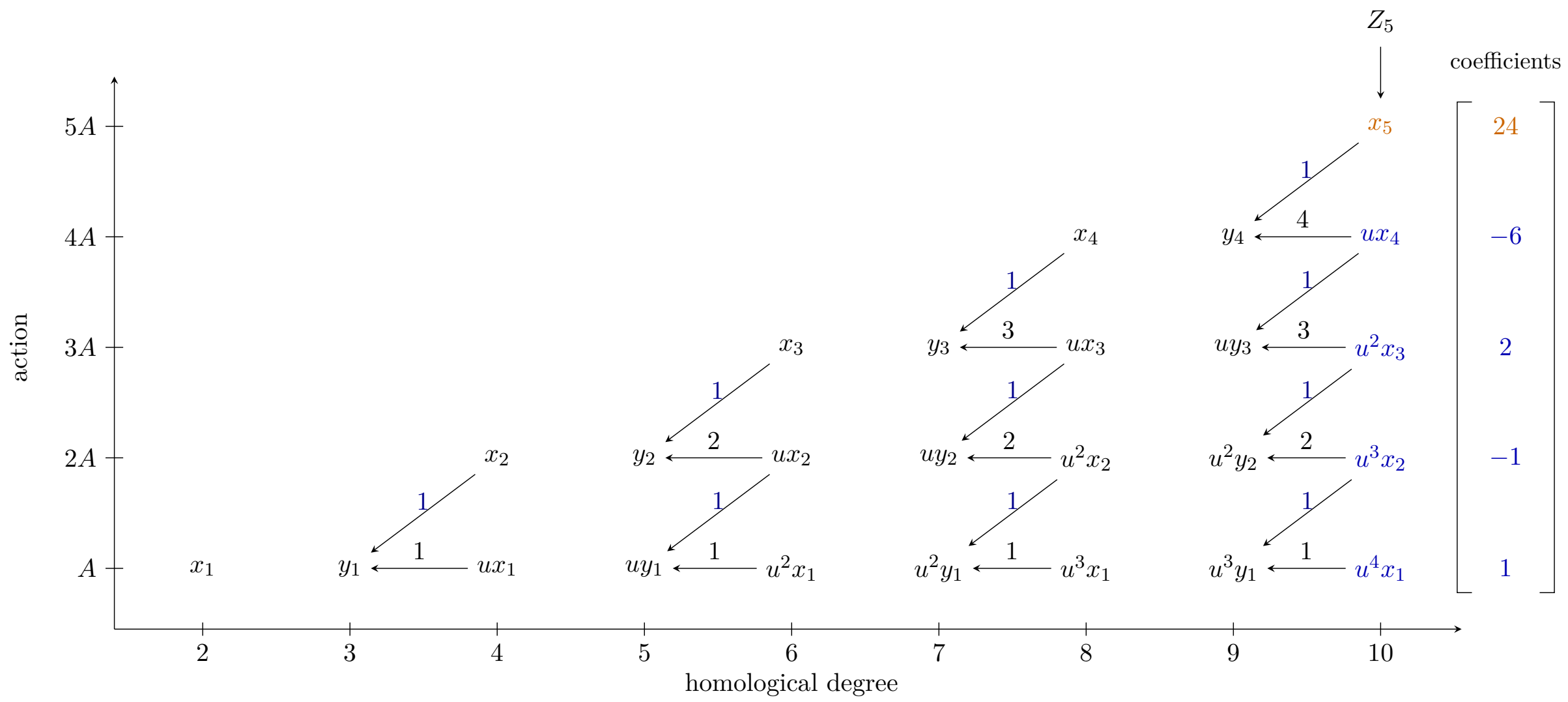

The stabilized obstruction is new as far as we know.
It is not detected by the capacities of McDuff--Siegel
\cite{mcduff2024symplectic} or their refinements due to Michala
\cite{michala2025refined}.
We explain the arithmetic behind the obstruction before describing
the sublevel calculation.  For an ellipsoid $E(a)$, the axial Reeb
orbits are the coordinate circles.  The orbit on the $i$-th axis has
action $a_i$, and its $r$-fold iterate has action $ra_i$ and covering
multiplicity $r$.  At action at most $L$, the number of these iterates
and the product of their covering multiplicities are
\[
 N_a(L)=\sum_i\left\lfloor\frac{L}{a_i}\right\rfloor,
 \qquad
 \Delta_a(L)=\prod_i\left\lfloor\frac{L}{a_i}\right\rfloor!.
\]
The rational capacity is the first level with at least $k$ axial
iterates \cite[Remark~2.2]{gutt2018symplectic}.  With cyclic
coefficients there is a second way for the distinguished class
to appear: the product of multiplicities can become zero in the
coefficient ring.  Define
\[
 \tau_{p^e}(a)=\inf\{L>0\mid p^e\mid\Delta_a(L)\}.
\]
Our ellipsoid formula, \cref{thm:cyclic-ellipsoid-capacities}, is
\[
 c_k^{\bbZ/p^e}(E(a))
       =\min\{c_k^{\bbQ}(E(a)),\tau_{p^e}(a)\}.
\]
This extends the prime-field calculation to all finite cyclic
coefficients.  In \cref{fig:disk-arithmetic-complex}, this distinction
can already be seen for the disk $D(A)=E(A)$.  With the coefficients
shown alongside the rightmost column, the blue terms have boundary
$-24y_4$: the horizontal arrows contribute the product
$1\cdot2\cdot3\cdot4=4!$.  Modulo $8$, these terms already form a
cycle at action $4A$, whereas over $\bbQ$ the orange term at action
$5A$ is needed to cancel the boundary.  Thus
\[
 c_5^{\bbZ/8}(D(A))=4A<5A=c_5^{\bbQ}(D(A)).
\]
At the blue cutoff there are only four iterates, and none of their
covering multiplicities is zero modulo $8$, but their product is.

% Positive equivariant disk complex from Zhao, Section 8.1.
% For covering multiplicity r, the lower generator has degree 2r.
% Horizontal coefficients are covering multiplicities; diagonal coefficients are 1.

For ellipsoids, \eqref{eq:ellipsoid-threshold-rational-sequence}
expresses the arithmetic threshold as a monotone function of the
rational capacity sequence.  Thus all cyclic capacity inequalities
between ellipsoids are implied by the rational ones; see
\cref{thm:no-new-cyclic-ellipsoid-obstruction}.
To obtain an additional obstruction we pass to more general convex
toric domains.  Their inclusions into supporting ellipsoids give
an upper bound with the same arithmetic contribution.

For a moment region $\Omega$, write
\[
 h_\Omega(v)=\sup_{w\in\Omega}\langle v,w\rangle,
 \qquad \bbN=\{0,1,2,\ldots\}.
\]
The rational toric formula of Gutt and Hutchings is
\[
 c_k^{\bbQ}(X_\Omega)
       =\min_{\substack{v\in\bbN^n\\\sum_i v_i=k}}h_\Omega(v).
\]
The corresponding prime-power threshold is
\[
 T_{p^e}(\Omega)
       =\min_{\substack{v\in\bbN^n\\p^e\mid\prod_i v_i!}}h_\Omega(v),
\]
and \cref{prop:cyclic-supporting-ellipsoid-bound} gives
\[
 c_k^{\bbZ/p^e}(X_\Omega)
       \leq\min\{c_k^{\bbQ}(X_\Omega),T_{p^e}(\Omega)\}.
\]
For $e<p+2$, these thresholds can themselves be expressed in terms
of rational capacities; see \cref{prop:small-exponent-threshold}.
The calculation for $\bbZ/243$ suggests a pattern for odd primes
$p$ with coefficients in $\bbZ/p^{p+2}$, discussed at the end of
\cref{sec:243-domains}.  We prove the case $p=3$, where the analytic
and topological mechanism are simpler; the other cases remain
conjectural.

The rational capacities minimize the support function over lattice
points of each fixed coordinate sum, while $T_{p^e}$ minimizes it
over the single set selected by the factorial divisibility condition.
These two minimizations can order domains differently.  For our
source and target, all rational capacity inequalities hold, whereas
\[
 T_{243}(\Omega_X)=\frac{48}{5}>9=T_{243}(\Omega_Y).
\]
For this pair, \cref{thm:243-obstruction} identifies these thresholds
with $c_{15}^{\bbZ/243}(X)$ and $c_{15}^{\bbZ/243}(Y)$, respectively,
so their reversed order gives the embedding obstruction.

We chose $X$ to obtain this reversal of the arithmetic thresholds
while keeping the topology needed for the capacity calculation tractable.
The supporting-ellipsoid bound does not show that the source
capacity attains this value.  A rational homology calculation can use cancellations
which no longer hold with torsion coefficients.  To determine the
image with the new coefficients, we use the topology of the filtered
spaces themselves.

The Clarke--Floer comparison of Abbondandolo and Kang
\cite{abbondandolo2022clarke} and its equivariant version due to
Matijevi\'c \cite{matijevic2024positive} make this possible.
For a convex body $C$ containing the origin in its interior, not
necessarily toric, let $\Lambda$ be the space of mean-zero $H^1$ loops
of period one and
symplectic action one.  Clarke's dual functional on this space is
\[
 \Psi_C(x)=\int_{\bbR/\bbZ}H_C^*(-J_0\dot x(t))\,dt,
 \qquad
 \Lambda_C^{<L}=\{x\in\Lambda\mid\Psi_C(x)<L\},
\]
where $H_C^*$ is the Legendre transform of the positively
two-homogeneous Hamiltonian equal to one on $\partial C$.
We give the definitions in \cref{sec:variational-comparison}.
Matijevi\'c works over $\bbF_2$ for simplicity and notes that his
identification of capacities holds over arbitrary fields
\cite[discussion following Theorem~A]{matijevic2024positive}.
For the cyclic coefficients considered here, we give the signed
version of his comparison over $\bbZ$ in
\cref{sec:coefficient-comparison}.  Passing to a commutative unital
coefficient ring $R$ then gives
\[
 H_*^{S^1}(\Lambda_C^{<L};R)
       \cong \SH_{*+n+1}^{S^1,+,<L}(C;R).
\]
Here $C$ is smooth, non-degenerate and strongly convex, and $L$
lies outside its action spectrum.  The compatibility with filtration
and limiting maps identifies the capacity with the first level at
which $\beta_{k-1}$ is in the image of the Borel projection
\[
 H_{2k-2}^{S^1}(\Lambda_C^{<L};R)
       \longrightarrow H_{2k-2}(BS^1;R).
\]
Convex approximation extends this characterization to the domains
with corners used in our examples.

For every ellipsoid, we give an $S^1$-equivariant homotopy equivalence
of filtered spaces:
\[
 \Lambda_{E(a)}^{<L}\simeq_{S^1}S(V_L),
 \qquad
 V_L=\bigoplus_{\substack{1\leq i\leq n,\ r\geq1\\ra_i<L}}\bbC_r.
\]
Here $\bbC_r$ is the line on which time shift by $\theta$ acts by
$e^{-2\pi ir\theta}$.
The associated bundle's Euler class has degree determined by the orbit count and
coefficient determined by the product of covering multiplicities,
which explains the two thresholds in the ellipsoid formula.

For the source $X$, the dual Hamiltonian is a pointwise minimum of
quadratic forms.  The minimizing parameter may vary along a loop,
preventing a direct Fourier calculation.  An equivariant homotopy
reparametrizes each coordinate at constant speed, preserving its
symplectic action and not increasing the dual functional.  This
allows a single parameter to realize the minimum along the whole
loop, replacing the sublevel by a union of quadratic sublevels.
Fourier damping then identifies the strict sublevel at $48/5$ with
\[
 S(V_{9,2})\cup_{S(V_{8,2})}S(V_{8,6}),
 \qquad
 V_{r,s}=\bigoplus_{j=1}^r\bbC_j\oplus\bigoplus_{j=1}^s\bbC_j.
\]
The Euler classes determine the maps between the associated Thom
spaces.  Their homotopy pushout computes the Borel image and shows
that $\beta_{14}$ has not yet appeared at this strict level.
Together with the upper bound, this proves
$c_{15}^{\bbZ/243}(X)=48/5$.
The maximal lattice rectangles in the positive polar suggest the
same two-sphere pattern for the odd-prime family.

To prove \cref{thm:capacity-disk-stabilization}, we work at the level
of the equivariant homotopy type.  Damping a loop's disk coordinate
gives an $S^1$-equivariant strong deformation retraction
\[
 \Lambda_{C\times D(A)}^{<L}\searrow\Lambda_C^{<L},
 \qquad 0<L<A.
\]
The retractions are compatible with the sublevel inclusions, so
the capacities are preserved for every coefficient ring.
Applying this to compact subsets of a stabilized embedding gives
the lower bound in \cref{cor:stabilized-243}; inclusion into
$E(16/15,16/5)$ attains the bound.

The original motivation for this work came from homotopical aspects
of symplectic invariants and led to these coefficient-sensitive
calculations.  The work of Kawamoto and Shelukhin
\cite{kawamoto2023spectral} provides another instance in which
integer coefficients in Hamiltonian Floer theory detect phenomena
not seen over fields.  The Clarke framework applies to arbitrary
convex domains, although explicit sublevel descriptions may require
further analytic ideas.

\paragraph{Organization of the paper.}
In \cref{sec:cyclic-capacities}, we recall the coefficient capacities
and state the ellipsoid formula and its consequence for embeddings
between ellipsoids.
In \cref{sec:cyclic-coefficient-capacities}, we derive the toric
upper bound, prove \cref{prop:prime-field-toric-capacities} and
describe the small-exponent thresholds.
In \cref{sec:variational-comparison}, we introduce the Clarke
functional and determine the ellipsoid sublevels.
The coefficient comparison in \cref{sec:coefficient-comparison}
identifies the capacities and completes the ellipsoid calculation.
In \cref{sec:243-domains}, we prove the torsion obstruction by the
two-sphere calculation and discuss the patterns it suggests for
other coefficients.
Finally, \cref{sec:stabilization-sketch} proves the general stabilization
theorem and applies it to obtain the sharp stabilized embedding threshold.

\paragraph{Use of AI.}
The proofs of
\cref{lem:243-reparametrization,lem:Clarke-disk-stabilization}
were developed with the assistance of AI. AI also helped identify
and clarify the relevant arguments in the literature used in the
coefficient comparison in \cref{prop:integral-matijevic-comparison}.
Further uses included manuscript revision and review, preparation
of figures, and locating references. The models used belong to the
GPT-6 family. The author takes full responsibility for the
mathematical content and its presentation.

\paragraph{Acknowledgments.}
I would like to thank Tomer Schlank and Yoel Groman for directing me
toward symplectic geometry in general and capacities in particular.
I thank Reut Medion and Hagit and Meir Frankel for their hospitality
while I was working on this research.
I thank Leonid Polterovich for wondering why not just look at ordinary
homology first.
The author was partially supported by the British Friends Graduate
Scholarship Fund.

\section{Cyclic capacities}
\label{sec:cyclic-capacities}

Changing the coefficient ring changes the action level at which a
distinguished homology class appears.  We study this dependence for
positive $S^1$-equivariant symplectic homology, focusing on
coefficients $\bbZ/p^e$.  The ellipsoid formula below exhibits the arithmetic
contribution of the orbit multiplicities.

\subsection{Definition}
\label{subsec:coefficient-capacity-definition}

Let $C\subset\bbC^n$ be a smooth convex domain containing the origin, and
let $R$ be a commutative unital ring. We briefly recall the filtered
homology used to define the capacities, following Gutt and Hutchings
\cite[Sections~3--5]{gutt2018symplectic}.
Hamiltonian Floer complexes on the completion of $C$ are generated by
periodic orbits and filtered by their action. Passing to the direct
limit over admissible Hamiltonians gives symplectic homology.
The $S^1$-equivariant version incorporates time shift of loops, and
its positive part is obtained by quotienting out the low-action
subcomplex. We write its action sublevels as
\[
    \SH_*^{S^1,+,<L}(C;R).
\]
Increasing $L$ gives maps between these groups and a map to their
direct limit $\SH_*^{S^1,+}(C;R)$.

The low-action part is relative homology of $(C,\partial C)$,
tensored with $H_*(BS^1;R)$, with the usual grading shift.
Since $C$ is contractible and displaceable, its full equivariant
symplectic homology vanishes. Thus
the connecting map for the low-action/positive exact sequence
identifies positive equivariant symplectic homology with a shifted
copy of $H_*(BS^1;R)$. We use the normalization
\[
 \SH_{*+n+1}^{S^1,+}(C;R)\cong H_*(BS^1;R).
\]
Let $\beta_j\in H_{2j}(BS^1;R)$ be dual to $u^j$, where $u$ is a
generator of $H^2(BS^1;\bbZ)$, mapped to $R$ coefficients.
Composing the filtered-to-limit map with this identification gives
\[
 \rho_{L,R}^C\colon
 \SH_{*+n+1}^{S^1,+,<L}(C;R)\longrightarrow H_*(BS^1;R).
\]
The coefficient-$R$ capacity is
\begin{equation}
 c_k^R(C)=\inf\{L>0\mid\beta_{k-1}\in\im\rho_{L,R}^C\},
 \qquad k\geq1.
 \label{eq:coefficient-capacity-definition}
\end{equation}
By exactness, this image condition also says that the low-action
class $[C,\partial C]\otimes\beta_{k-1}$ becomes zero in full
equivariant symplectic homology at that action level.
Thus $\beta_{k-1}$ corresponds to degree $2k+n-1$ in positive
equivariant symplectic homology.  With rational coefficients these are
the Gutt--Hutchings capacities \cite{gutt2018symplectic}.
Their transfer argument \cite[Lemma~4.2(a),(c)]{gutt2018symplectic},
with the coefficient conventions used here, gives monotonicity and
conformality.  We extend the capacities to convex bodies by the
following standard argument.

\paragraph{Convex approximation.}
We recall the approximation argument used below to extend capacity
identities and inequalities to domains with corners.  All convex bodies
in this discussion are compact and contain the origin in their interiors.

Let $c$ be a finite nonnegative invariant of smooth non-degenerate
strongly convex domains, monotone under inclusions and satisfying
\[
 c(rC)=r^2c(C),\qquad r>0.
\]
For an arbitrary convex body $C$, extend $c$ by taking the supremum
over smooth non-degenerate strongly convex bodies inside it:
\[
 c(C)=\sup_{K\subset\operatorname{int}(C)}c(K).
\]
This agrees with the original invariant on smooth non-degenerate
strongly convex domains: the upper bound follows from monotonicity,
and the lower bound follows by taking $K=rC$ and letting $r\uparrow1$.
The extension retains inclusion monotonicity and conformality.

For every $0<\varepsilon<1$, choose a smooth non-degenerate strongly
convex approximation satisfying
\[
 (1-\varepsilon)C\subset C_\varepsilon\subset(1+\varepsilon)C.
\]
Monotonicity and conformality then give
\[
 (1-\varepsilon)^2c(C)
 \leq c(C_\varepsilon)
 \leq(1+\varepsilon)^2c(C).
\]
Hence $c(C_\varepsilon)\to c(C)$ as $\varepsilon\to0$.
The same estimate applies to any invariant already defined on convex
bodies with these two properties.  Thus agreement on smooth
non-degenerate strongly convex domains implies agreement on all convex
bodies.  Inequalities involving finitely many capacity values and
continuous operations, such as minima and maxima, also pass to the limit.

Symplectic monotonicity is preserved as well.  Suppose that $c$ is
monotone under symplectic embeddings of smooth domains and that
\[
 \varphi\colon\operatorname{int}(C)\hookrightarrow\operatorname{int}(D)
\]
is symplectic.  Choose smooth non-degenerate strongly convex bodies
$C_\varepsilon^-$ and $D_\varepsilon^+$ with
\[
 \begin{aligned}
 (1-\varepsilon)C
 &\subset C_\varepsilon^-\subset\operatorname{int}(C),\\
 D&\subset\operatorname{int}(D_\varepsilon^+)
       \subset D_\varepsilon^+\subset(1+\varepsilon)D.
 \end{aligned}
\]
The restriction of $\varphi$ embeds $C_\varepsilon^-$ into
$\operatorname{int}(D_\varepsilon^+)$.  It is exact, as required for
the transfer map, since $\varphi^*\lambda_0-\lambda_0$ is closed and
$C_\varepsilon^-$ is contractible.  Hence
\[
 (1-\varepsilon)^2c(C)
 \leq c(C_\varepsilon^-)
 \leq c(D_\varepsilon^+)
 \leq(1+\varepsilon)^2c(D).
\]
Letting $\varepsilon\to0$ gives $c(C)\leq c(D)$.

For finite cyclic coefficients, it suffices to work with prime powers.
Indeed, the Chinese remainder decomposition
$\bbZ/m\cong\prod_{p^e\parallel m}\bbZ/p^e$ splits the filtered
complex and the distinguished class componentwise.  The class appears
precisely when all its components appear, so
\begin{equation}
 \label{eq:cyclic-CRT-capacity}
 c_k^{\bbZ/m}(C)=\max_{p^e\parallel m}c_k^{\bbZ/p^e}(C).
\end{equation}
Convex approximation preserves this identity.
We therefore take $p$ prime and $e\geq1$ in the calculations below.

\subsection{Ellipsoids}

For $a=(a_1,\ldots,a_n)$ with $a_i>0$, write
\[
 E(a)=\left\{z\in\bbC^n\ \middle|\
              \sum_{i=1}^n\frac{\pi|z_i|^2}{a_i}\leq1\right\}.
\]
The axial Reeb orbits are the coordinate circles on $\partial E(a)$.
The orbit on the $i$th coordinate axis has action $a_i$; its $r$-fold
iterate has action $ra_i$ and covering multiplicity $r$.

For $L>0$, the number of axial iterates with action at most $L$
and the product of their covering multiplicities are
\[
 N_a(L)=\sum_{i=1}^n\left\lfloor\frac{L}{a_i}\right\rfloor,
 \qquad
 \Delta_a(L)=\prod_{i=1}^n\left\lfloor\frac{L}{a_i}\right\rfloor!.
\]
Here $0!=1$. Gutt and Hutchings
\cite[Remark~2.2]{gutt2018symplectic} give the rational capacity as
\begin{equation}
 c_k^{\bbQ}(E(a))=\min\{L>0\mid k \leq N_a(L)\}.
 \label{eq:rational-ellipsoid-count}
\end{equation}
The arithmetic threshold records when the product of covering
multiplicities becomes divisible by $p^e$:
\begin{equation}
 \tau_{p^e}(a)=\inf\{L>0\ \mid p^e\mid\Delta_a(L)\}.
 \label{eq:cyclic-ellipsoid-threshold}
\end{equation}
For a prime $p$, we write $\nu_p$ for the $p$-adic valuation.

\begin{theorem}
 \label{thm:cyclic-ellipsoid-capacities}
 For every ellipsoid $E(a)$, every prime $p$, and all $e,k\geq1$,
 \begin{equation}
  c_k^{\bbZ/p^e}(E(a))=\min\{c_k^{\bbQ}(E(a)),\tau_{p^e}(a)\}.
  \label{eq:cyclic-ellipsoid-formula}
 \end{equation}
\end{theorem}

The divisibility condition can equivalently be written using valuations:
\begin{equation}
 \tau_{p^e}(a)=
 \inf\left\{L>0\ \middle|\
       \sum_{i=1}^n\nu_p\left(
           \left\lfloor\frac{L}{a_i}\right\rfloor!\right)\geq e\right\}.
 \label{eq:prime-power-ellipsoid-threshold}
\end{equation}

The cyclic sequence agrees with the rational sequence until it reaches
$\tau_{p^e}(a)$, and is constant thereafter.

For a moment region $\Omega$, its positive polar is
\[
 \Omega^*=\{v\in\bbR_{\geq0}^n\mid h_\Omega(v)\leq1\}.
\]
The moment simplex $\Omega$ of $E(a)$ has positive polar
$\prod_i[0,1/a_i]$. The largest lattice point in $L\Omega^*$
has coordinates $v_i=\lfloor L/a_i\rfloor$, so
$N_a(L)=\sum_i v_i$ and $\Delta_a(L)=\prod_i v_i!$.
The rational capacity is the first scale at which $\sum_i v_i\geq k$;
with coefficients in $\bbZ/p^e$, the condition $p^e\mid\prod_i v_i!$
also contributes.

The formula can also be obtained directly from the positive
$S^1$-equivariant chain complex.  We use the Clarke model here to
explain how covering multiplicities enter the homological image and to
prepare for the two-sphere calculation, completing the proof in
\cref{sec:coefficient-comparison} after establishing the comparison
over coefficient rings.
We first record its consequence for embeddings between ellipsoids.

\begin{theorem}
    \label{thm:no-new-cyclic-ellipsoid-obstruction}

    Let $E(a)$ and $E(b)$ be ellipsoids.  If
    \[
        c_k^{\bbQ}(E(a))
        \leq
        c_k^{\bbQ}(E(b))
        \qquad\text{for every }k,
    \]
    then, for every prime $p$ and all $e,k\geq1$,
    \[
        c_k^{\bbZ/p^e}(E(a))
        \leq
        c_k^{\bbZ/p^e}(E(b)).
    \]
\end{theorem}

\begin{proof}
    By \eqref{eq:rational-ellipsoid-count}, $N_a(L)$ is the number of
    indices $k$ with $c_k^{\bbQ}(E(a))\leq L$. The hypothesis implies
    that each index counted for $E(b)$ is also counted for $E(a)$, so
    \[
        N_a(L)\geq N_b(L)\qquad(L>0).
    \]
    Legendre's formula therefore yields, for every prime $p$,
    \[
        \begin{aligned}
            \nu_p(\Delta_a(L))
            &=\sum_i\nu_p\left(\left\lfloor\frac{L}{a_i}\right\rfloor!\right)\\
            &=\sum_{s\geq1}\sum_i
              \left\lfloor\frac{L}{p^s a_i}\right\rfloor
            =\sum_{s\geq1}N_a\left(\frac{L}{p^s}\right)
            \geq
            \sum_{s\geq1}N_b\left(\frac{L}{p^s}\right)
            &=\sum_{s\geq1}\sum_i
              \left\lfloor\frac{L}{p^s b_i}\right\rfloor\\
            &=\sum_i\nu_p\left(\left\lfloor\frac{L}{b_i}\right\rfloor!\right)
            =\nu_p(\Delta_b(L)).
        \end{aligned}
    \]
    Thus $\Delta_b(L)\mid \Delta_a(L)$, and consequently
    $\tau_{p^e}(a)\leq\tau_{p^e}(b)$.  The result now follows from
    \eqref{eq:cyclic-ellipsoid-formula}.
\end{proof}

By \eqref{eq:cyclic-CRT-capacity}, the same conclusion holds for every
finite cyclic coefficient ring.

\begin{remark}
The identity above expresses $\nu_p(\Delta_a(L))$ as the number
of entries at most $L$ in the list below.  Thus $\tau_{p^e}(a)$
is its $e$-th entry, counted with multiplicity.
\begin{equation}
 \label{eq:ellipsoid-threshold-rational-sequence}
 \bigsqcup_{s\geq1}
 \bigl\{p^s c_j^{\bbQ}(E(a)):j\geq1\bigr\},
 \qquad\text{arranged in nondecreasing order.}
\end{equation}
This description makes its monotonicity in the rational capacity
sequence explicit and, together with \eqref{eq:cyclic-ellipsoid-formula},
explains why cyclic coefficients give no additional embedding
obstructions between ellipsoids.
\end{remark}

\section{From ellipsoids to convex toric domains}
\label{sec:cyclic-coefficient-capacities}

A convex toric domain lies in each of its supporting ellipsoids.
The ellipsoid formula therefore gives an arithmetic upper bound
for its cyclic capacities.  Over a prime field, the low-action
comparison of \c{C}ineli, Ginzburg, and G\"urel makes this bound exact,
proving \cref{prop:prime-field-toric-capacities}.  The resulting
formula shows why changing to prime-field coefficients gives no
additional embedding obstructions between convex toric domains.
We then examine the thresholds for higher prime powers, before
turning to the sublevel calculation needed to prove our obstruction.

\subsection{The supporting-ellipsoid bound}

Write $X_\Omega=\mu^{-1}(\Omega)$, where
$\mu(z)=(\pi|z_1|^2,\ldots,\pi|z_n|^2)$.
We use the convention for convex toric domains of Gutt and Hutchings
\cite{gutt2018symplectic}: $\Omega\subset\bbR_{\geq0}^n$ is compact
with nonempty interior, and its reflection across all coordinate
hyperplanes is convex.  This condition makes $\Omega$ convex and
downward closed in each coordinate, so $X_\Omega$ is a compact
convex body in $\bbC^n$ containing the origin in its interior.
With $\bbN=\{0,1,2,\ldots\}$, put
\[
    h_\Omega(v)=\max_{w\in\Omega}\langle v,w\rangle.
\]
For convex toric domains, the rational Gutt--Hutchings formula
\cite{gutt2018symplectic} is
\[
    c_k^{\bbQ}(X_\Omega)
    =\min_{\substack{v\in\bbN^n\\\sum_i v_i=k}}h_\Omega(v).
\]

For a convex toric domain $X_\Omega\subset\Cn$, define
\begin{equation}
    T_{p^e}(\Omega)
    =
    \min_{\substack{v\in\bbN^n\\
                    p^e\mid\prod_i v_i!}}
    h_\Omega(v).
    \label{eq:def-cyclic-toric-threshold}
\end{equation}

\begin{proposition}
    \label{prop:cyclic-supporting-ellipsoid-bound}
    For every convex toric domain $X_\Omega$, every prime $p$, and all $e,k\geq1$,
    \begin{equation}
        c_k^{\bbZ/p^e}(X_\Omega)
        \leq
        \min\{c_k^{\bbQ}(X_\Omega),T_{p^e}(\Omega)\}.
        \label{eq:cyclic-toric-upper-bound}
    \end{equation}
\end{proposition}

The bound can also be written as
    \begin{equation}
        c_k^{\bbZ/p^e}(X_\Omega)
        \leq
        \min_{\substack{v\in\bbN^n\\
          \sum_i v_i\geq k\ \text{or}\ 
          p^e\mid\prod_i v_i!}}
        h_\Omega(v).
        \label{eq:cyclic-single-support-bound}
    \end{equation}

The two alternatives record the orbit count and the divisibility of the
product of covering multiplicities, respectively.

\begin{proof}
    Fix $v\in\bbN^n$ satisfying either condition in
    \eqref{eq:cyclic-single-support-bound}.  First suppose that every
    $v_i>0$, and put
    \[
        L=h_\Omega(v),\qquad a_i=\frac{L}{v_i}.
    \]
    The support-function inequality gives
    \[
        \sum_i\frac{w_i}{a_i}
        =\frac{\langle v,w\rangle}{L}\leq1
        \qquad(w\in\Omega),
    \]
    so $X_\Omega\subset E(a)$. At action $L$, this ellipsoid has
    exactly $v_i$ axial iterates from the $i$th axis, so
    \[
        N_a(L)=\sum_i v_i,\qquad \Delta_a(L)=\prod_i v_i!.
    \]
    Under either condition in \eqref{eq:cyclic-single-support-bound},
    the ellipsoid formula gives $c_k^{\bbZ/p^e}(E(a))\leq L$.
    Monotonicity therefore yields
    \[
        c_k^{\bbZ/p^e}(X_\Omega)
        \leq c_k^{\bbZ/p^e}(E(a))\leq h_\Omega(v).
    \]

    If some $v_i=0$, use $v^\varepsilon=v+\varepsilon(1,\ldots,1)$
    with $0<\varepsilon<1$ in the supporting-ellipsoid construction.
    At action $h_\Omega(v^\varepsilon)$, the corresponding ellipsoid
    still has exactly $v_i$ axial iterates from the $i$th axis, since
    $\lfloor v_i+\varepsilon\rfloor=v_i$.  The same argument therefore gives
    \[
        c_k^{\bbZ/p^e}(X_\Omega)
        \leq h_\Omega(v^\varepsilon)\longrightarrow h_\Omega(v)
        \qquad(\varepsilon\downarrow0).
    \]
    Minimizing over $v$ proves \eqref{eq:cyclic-single-support-bound};
    separating its two alternatives gives \eqref{eq:cyclic-toric-upper-bound}.
\end{proof}

The two quantities in \eqref{eq:cyclic-toric-upper-bound} measure
different properties of the same support function.  The rational
capacity uses the coordinate sum, while $T_{p^e}$ uses the divisibility
of the product of factorials.

\subsection{Prime-field capacities and prime-power thresholds}

We now prove \cref{prop:prime-field-toric-capacities}.  For prime
coefficients the supporting-ellipsoid bound is exact, and its
arithmetic threshold is determined by the first rational capacity.

\begin{proof}[Proof of \cref{prop:prime-field-toric-capacities}]
The condition $p\mid\prod_i v_i!$ is equivalent to $v_i\geq p$
for some $i$.  By coordinatewise monotonicity of $h_\Omega$,
\[
 T_p(\Omega)=p\min_i h_\Omega(e_i)
             =p\,c_1^{\bbQ}(X_\Omega),
\]
where $e_i$ are the standard coordinate vectors.
Thus \cref{prop:cyclic-supporting-ellipsoid-bound} gives the upper bound.

For the reverse inequality we use the coefficient comparison in
\cite[Section~2.4, discussion following~(2.21)]{cineli2026filtered}.
For a smooth non-degenerate strongly convex domain $C$, it gives
\[
 c_k^{\bbF_p}(C)<p\,a_{\min}(C)
       \quad\Longrightarrow\quad
 c_k^{\bbF_p}(C)=c_k^{\bbQ}(C),
\]
where $a_{\min}(C)$ is the least action of a closed Reeb orbit.
For convex domains this action equals $c_1^{\bbQ}(C)$
\cite{matijevic2024positive}.  Hence either
$c_k^{\bbF_p}(C)\geq p\,c_1^{\bbQ}(C)$ or it equals the rational
capacity.  In both cases,
\[
 c_k^{\bbF_p}(C)\geq
       \min\{c_k^{\bbQ}(C),p\,c_1^{\bbQ}(C)\}.
\]
Convex approximation extends this inequality to $X_\Omega$,
proving \eqref{eq:prime-field-toric-capacities}.
If all rational capacity inequalities hold for $X$ and $Y$, then
both entries in the minimum for $X$ are at most the corresponding
entries for $Y$.  This proves the asserted prime-field inequalities.
\end{proof}

This also answers, for convex toric domains, the calculation and
stabilization questions raised in
\cite[Section~1.2.1]{cineli2026filtered}.  The sequence is the rational
sequence truncated at $p\,c_1^{\bbQ}$, and eventually equals this value.
For a smooth domain this is the least action of a $p$-fold iterate.
Its first strict difference from the rational sequence occurs at
the first index with $c_k^{\bbQ}>p\,c_1^{\bbQ}$.
The formula also shows that, for each fixed $k$ and domain, the
prime-field capacity is nondecreasing as a function of the prime $p$.
These conclusions concern convex toric domains; the questions in
that paper are posed for a broader class of domains.

For higher prime powers, Legendre's formula writes the threshold as
\begin{equation}
 \label{eq:def-prime-power-toric-threshold}
 T_{p^e}(\Omega)=
 \min_{\substack{v\in\bbN^n\\\sum_i\nu_p(v_i!)\geq e}}h_\Omega(v).
\end{equation}
For small exponents, it is still determined by the rational capacities.

\begin{proposition}
\label{prop:small-exponent-threshold}
For every convex toric domain $X_\Omega$ and every prime $p$,
\begin{equation}
 \label{eq:small-exponent-threshold}
 \begin{aligned}
 T_{p^e}(\Omega)&=p\,c_e^{\bbQ}(X_\Omega)
                           &&(1\leq e\leq p),\\
 T_{p^{p+1}}(\Omega)
   &=\min\{p\,c_{p+1}^{\bbQ}(X_\Omega),
                        p^2c_1^{\bbQ}(X_\Omega)\}.
 \end{aligned}
\end{equation}
\end{proposition}

\begin{proof}
If $w\in\bbN^n$ has coordinate sum $e$, then $v=pw$ satisfies
$\sum_i\nu_p(v_i!)\geq e$.  Choosing $w$ to minimize the rational
support formula gives $T_{p^e}\leq p\,c_e^{\bbQ}$.
Also, a coordinate vector $p^2e_i$ has factorial valuation $p+1$,
giving $T_{p^{p+1}}\leq p^2c_1^{\bbQ}$.

For the reverse inequalities, let $v$ be admissible for $T_{p^e}$.
If all coordinates are less than $p^2$, Legendre's formula reduces to
\[
 \sum_i\nu_p(v_i!)=\sum_i\left\lfloor\frac{v_i}{p}\right\rfloor.
\]
There is then $w\in\bbN^n$ with $\sum_iw_i=e$ and $pw\leq v$
coordinatewise, so
$h_\Omega(v)\geq p h_\Omega(w)\geq p\,c_e^{\bbQ}$.
If some coordinate is at least $p^2$, then
$h_\Omega(v)\geq p^2c_1^{\bbQ}$.
For $e\leq p$, we also have
$c_e^{\bbQ}\leq e\,c_1^{\bbQ}\leq p\,c_1^{\bbQ}$,
by using $e$ times a coordinate vector in the rational formula.
Thus this second case gives $h_\Omega(v)\geq p\,c_e^{\bbQ}$
when $e\leq p$, and the other alternative in the asserted minimum
when $e=p+1$.
\end{proof}

Our example uses $p=3$ and $e=5$, the next exponent $e=p+2$.
Here the distinction between coordinate sums and factorial valuations leads to a
different ordering of source and target thresholds.

To turn the arithmetic upper bound into an obstruction, we must
determine the filtered homological image for the source.
The cancellations available with rational coefficients need not
persist over $\bbZ/p^e$: an integer which is invertible over $\bbQ$
may be a nonunit in the new coefficient ring.
We therefore pass to the Clarke sublevels themselves.  Their
equivariant topology retains the covering multiplicities and allows
us to prove that the upper bound is attained for the source below.

\section{Clarke duality and ellipsoid sublevels}
\label{sec:variational-comparison}

Clarke's dual functional \cite{clarke1979classical} gives access to
the topology of the filtered spaces themselves.
Matijevi\'c identifies their equivariant homology with positive
equivariant symplectic homology over $\bbZ/2$
\cite{matijevic2024positive}, extending the non-equivariant comparison
of Abbondandolo--Kang \cite{abbondandolo2022clarke}.
The coefficient extension needed here is established in
\cref{sec:coefficient-comparison}.

After recalling the functional, we describe its ellipsoid sublevels as
spheres in finite-dimensional circle representations.  The classical
variational calculation of ellipsoid capacities appears in
Ekeland--Hofer \cite[Proposition~4]{ekeland1990symplecticII}.
The homological image supplies the
capacity formula, and the same representation spheres will be the
pieces of the source sublevel in our obstruction.

\subsection{The functional and its sublevels}

Let $\mathbb T=\bbR/\bbZ$ and identify $(\bbR^{2n},J_0)$ with $\bbC^n$.
On the mean-zero Sobolev loop space
\[
 H^1_0=\left\{x\in H^1(\mathbb T,\bbR^{2n})\ \middle|\
                   \int_{\mathbb T}x(t)\,dt=0\right\},
\]
we define the symplectic action by
\[
 \cA \colon H_0^1 \to \bbR, \qquad \cA(x)=\frac12\int_{\mathbb T}
              \langle J_0x(t),\dot x(t)\rangle\,dt,
 \qquad \Lambda=\cA^{-1}(1).
\]
The circle acts by $\theta\cdot x=x(\,\cdot-\theta)$.
For a convex body $C$ containing the origin in its interior, its gauge is
\[
 j_C(z)=\inf\{r>0\mid z\in rC\}.
\]
We take $H_C=j_C^2$, so $H_C$ is positively two-homogeneous and equals
one on $\partial C$.  Its Legendre transform is defined by
\[
 H_C^*(p)=\sup_{z\in\bbR^{2n}}
             \bigl\{\langle p,z\rangle-H_C(z)\bigr\}.
\]

The standard support-function identity
\cite[Section~2.3]{haimkislev2019symplectic} gives
\begin{equation}
 H_C^*(p)=\frac{s_C(p)^2}{4},
 \qquad s_C(p)=\sup_{u\in C}\langle p,u\rangle.
 \label{eq:Clarke-support-dual}
\end{equation}
We now define the dual energy
\begin{equation}
 \mathcal H_C(x)=\int_{\mathbb T}H_C^*(-J_0\dot x(t))\,dt,
 \qquad \Psi_C=\mathcal H_C|_\Lambda.
 \label{eq:Clarke-dual-functional}
\end{equation}
For smooth strictly convex $C$, critical points of $\Psi_C$ encode
closed characteristics on $\partial C$, with critical values equal
to their actions after normalization.  The same functional is defined
for convex domains with corners.

Write
\[
 \Lambda_C^{<L}=\{x\in\Lambda\mid\Psi_C(x)<L\}.
\]
The Borel projection is
\[
 p_L^C\colon(\Lambda_C^{<L})_{hS^1}\longrightarrow BS^1.
\]
In this model the corresponding spectral number is
\begin{equation}
 c_k^{R,\Psi}(C)=
 \inf\left\{L>0\ \middle|\
  \beta_{k-1}\in
  \im\bigl((p_L^C)_*\colon H_{2k-2}^{S^1}(\Lambda_C^{<L};R)
                    \longrightarrow H_{2k-2}(BS^1;R)\bigr)\right\}.
 \label{eq:Clarke-capacity-definition}
\end{equation}
The equality $c_k^{R,\Psi}(C)=c_k^R(C)$ is proved in
\cref{cor:Clarke-Floer-capacities} by identifying the Borel projection
with the map defining the Floer capacity.

For $m\geq1$, let $\bbC_m$ denote the complex line with action
$\theta\cdot z=e^{-2\pi im\theta}z$, and write
$\cL_m=ES^1\times_{S^1}\bbC_m$.  We choose the generator above as
$u=c_1(\cL_1)$, so $c_1(\cL_m)=mu$.
For a representation $V$, $S(V)$ denotes its unit sphere, with
$S(0)=\varnothing$.

\subsection{Ellipsoid sublevels}
\label{sec:ellipsoid-Clarke-homotopy}

For $E=E(a)$, rescaling the coordinates gives
\[
 s_E(p)
 =\sup_{\sum_j\pi|z_j|^2/a_j\leq1}\langle p,z\rangle
 =\left(\sum_{j=1}^n\frac{a_j|p_j|^2}{\pi}\right)^{1/2}.
\]
The last equality is Cauchy--Schwarz on the rescaled unit ball.
Thus \eqref{eq:Clarke-support-dual} gives
\[
 H_E^*(p)=\frac{s_E(p)^2}{4}
        =\sum_{j=1}^n\frac{a_j|p_j|^2}{4\pi}.
\]
The sublevel at $L$ retracts onto a sphere in the positive Fourier modes
with action below $L$.  As $L$ passes action values, these spheres grow
by adding the corresponding weight lines.

Throughout this subsection, let $L>0$ lie outside the action spectrum
$\sigma(\partial E)=\{ma_j\mid m\geq1,\ 1\leq j\leq n\}$. Set
\[
 V_L:=\bigoplus_{\substack{1\leq j\leq n,\ m\geq1\\ma_j<L}}\bbC_m, \quad
 \dim_{\bbC}V_L=N_a(L).
\]
\begin{theorem}
\label{thm:ellipsoid-weighted-sphere}
There is an $S^1$-equivariant homotopy equivalence of filtered spaces
\begin{equation}
 \label{eq:ellipsoid-weighted-sphere}
 \Lambda_E^{<L}\simeq_{S^1}S(V_L).
\end{equation}
\end{theorem}

\begin{proof}
We construct a strong deformation retraction onto $\Lambda\cap V_L$
and identify this space with $S(V_L)$ by diagonal rescaling.
Write a zero-mean loop in Fourier coordinates as
\[
 x_j(t)=\sum_{m\in\bbZ\setminus\{0\}}
 e^{2\pi imt}x_{j,m}.
\]
Parseval's identity gives
\begin{align}
 \cA(x)
 &=\frac12\int_{\mathbb T}\langle J_0x,\dot x\rangle\,dt
   =\pi\sum_{j=1}^n\sum_{m\neq0}m|x_{j,m}|^2,
 \label{eq:ellipsoid-action-Fourier}\\
 \mathcal H_E(x)&=
 \pi\sum_{j=1}^n\sum_{m\neq0}a_jm^2|x_{j,m}|^2.
 \label{eq:ellipsoid-H-Fourier}
\end{align}
Indeed, on the $m$-th mode one has
$\dot x_j=2\pi mJ_0x_j$ and
$-J_0\dot x_j=2\pi m x_j$.  Time shift multiplies $x_{j,m}$ by
$e^{-2\pi im\theta}$, so the mode is the representation $\bbC_m$.

Consider the quadratic form
\[
 F_L:=\mathcal H_E-L\cA
 =\pi\sum_{j=1}^n\sum_{m\neq0}
 (a_jm^2-Lm)|x_{j,m}|^2.
\]
The fact that
\[
      a_jm^2-Lm<0 \iff 0<ma_j<L
\]
implies that the negative space is precisely $V_L$.  If $W_L$ is the complementary
sum of Fourier modes, then
\[
 H^1_0=V_L\oplus W_L,
 \qquad F_L|_{V_L}<0,
 \qquad F_L|_{W_L}\geq0.
\]

Writing $x\in\Lambda_E^{<L}=\{x\in H^1_0\mid \cA(x)=1,\ F_L(x)<0\}$
as $x=v+w$ with $v\in V_L$ and $w\in W_L$, the condition $F_L(x)<0$
forces $v\neq0$.
For $s\in[0,1]$, let $y_s=v+(1-s)w$.  Orthogonality of the Fourier
modes gives
\[
 F_L(y_s)=F_L(v)+(1-s)^2F_L(w)\leq F_L(x)<0.
\]
As $\mathcal H_E$ is positive definite and $L>0$, this inequality implies
$\cA(y_s)>0$.  Thus
\[
 R_s(x):=\frac{y_s}{\sqrt{\cA(y_s)}}
\]
is well-defined and lies in $\Lambda_E^{<L}$.  The homotopy is equivariant,
$R_0=\mathrm{id}$, $R_1(\Lambda_E^{<L})\subset\Lambda\cap V_L$, and it fixes
$\Lambda\cap V_L$ pointwise.  It is therefore the required strong
deformation retraction.

On $V_L$, formula~\eqref{eq:ellipsoid-action-Fourier} identifies
$\Lambda\cap V_L$ with the ellipsoid
\[
 \left\{v\in V_L\ \middle|\
 \pi\sum_{\substack{1\leq j\leq n,\ m\geq1\\ma_j<L}}m|v_{j,m}|^2=1\right\}.
\]
The map $(v_{j,m})\mapsto(\sqrt{\pi m}\,v_{j,m})$ is an
$S^1$-equivariant diffeomorphism from this ellipsoid to $S(V_L)$.
These rescalings commute with $V_L\subset V_{L'}$.  Applying the two
deformation retractions to a filtration map shows that, up to equivariant
homotopy, it is the standard equatorial inclusion.
\end{proof}

The equivariant sphere description now lets us read off the
homological image for any coefficient ring.

The representation $V_L$ defines a complex bundle
\begin{equation}
 \xi_L=ES^1\times_{S^1}V_L
      \cong\bigoplus_{\substack{1\leq j\leq n,\ m\geq1\\ma_j<L}}\cL_m
      \longrightarrow BS^1.
 \label{eq:ellipsoid-xi-L}
\end{equation}
Its sphere bundle is $S(V_L)_{hS^1}$.  The resulting image calculation
will give the ellipsoid capacity formula
\eqref{eq:cyclic-ellipsoid-formula}.

\begin{proposition}
\label{prop:ellipsoid-Borel-sphere}
For every commutative unital ring $R$ and every degree $q$,
\begin{equation}
 \im(p_L^E)_*
 =\ker\left(e_R(\xi_L)\frown(-)\colon
       H_q(BS^1;R)\longrightarrow H_{q-2N_a(L)}(BS^1;R)\right),
 \label{eq:ellipsoid-homological-Euler-kernel}
\end{equation}
where
\begin{equation}
 e_R(\xi_L)=\Delta_a(L)\,u^{N_a(L)}.
 \label{eq:ellipsoid-Euler-class}
\end{equation}
\end{proposition}

\begin{proof}
The sphere equivalence in \cref{thm:ellipsoid-weighted-sphere} gives
$(\Lambda_E^{<L})_{hS^1}\simeq S(\xi_L)$ over $BS^1$.
The disk/sphere-bundle pair gives the cofiber sequence
\begin{equation}
 S(\xi_L)_+\longrightarrow BS^1_+
       \xrightarrow{z_L}\operatorname{Th}(\xi_L).
 \label{eq:ellipsoid-Thom-cofiber}
\end{equation}
Under the Thom isomorphism, $(z_L)_*$ is cap product with the Euler
class; see \cite[Section~3.2]{hatcherVectorBundles}.  Exactness gives
the kernel formula.  Since $c_1(\cL_m)=mu$, the Euler class is the
product of the weights times $u^{N_a(L)}$.  The weights on the $j$th
coordinate are $1,\ldots,\lfloor L/a_j\rfloor$, whose product is
$\lfloor L/a_j\rfloor!$, giving $\Delta_a(L)$.
\end{proof}

The orbit count determines the degree of the Euler class, and the
product of covering multiplicities determines its coefficient.
The sphere bundle, rather than the ordinary orbit space, retains these
multiplicities. Over $\bbZ/p^e$ their product can vanish, making the
image larger than it is over $\bbQ$.

\begin{remark}
The same proofs apply at spectral levels, provided the orbit count and
multiplicity product include only $ma_j<L$.  We use non-spectral levels
to keep the formulas simple; they suffice for the capacity calculations.
\end{remark}

\begin{remark}
In our grading, the low-action/positive exact sequence of
Gutt--Hutchings \cite{gutt2018symplectic} reads
\[
 \cdots\longrightarrow
 \SH_{q+n+1}^{S^1,+,<L}(E;R)\longrightarrow
 \SH_{q+n}^{S^1,\leq\varepsilon}(E;R)\longrightarrow
 \SH_{q+n}^{S^1,<L}(E;R)\longrightarrow\cdots,
\]
where $0<\varepsilon<\min\{L,a_1,\ldots,a_n\}$.  Its Clarke analogue is
the exact sequence of the pair
\[
 \cdots\longrightarrow
 H_q^{S^1}(\Lambda_E^{<L};R)\longrightarrow
 H_q^{S^1}(\Lambda;R)\longrightarrow
 H_q^{S^1}(\Lambda,\Lambda_E^{<L};R)\longrightarrow\cdots.
\]
Thus the sublevel, ambient space and relative pair play the roles of
positive, low-action and full filtered homology, respectively.
Since $\Lambda$ is contractible, the middle term is $H_q(BS^1;R)$;
for the ellipsoid, the relative term is
$\widetilde H_q(\operatorname{Th}(\xi_L);R)$.
\end{remark}

\section{The coefficient comparison}
\label{sec:coefficient-comparison}

Matijevi\'c formulates the filtered equivariant Clarke--Floer comparison
over $\bbZ/2$ \cite[Theorem~C.2]{matijevic2024positive}.
The underlying construction of Abbondandolo--Kang
\cite{abbondandolo2022clarke} combines the hybrid Morse--Floer approach
of Abbondandolo--Schwarz \cite{abbondandolo2006cotangent} with a
finite-dimensional reduction following Viterbo
\cite{viterbo1989equivariant}.

To use the sublevel calculations with torsion coefficients, we give
the signed version of this construction over $\bbZ$ and then pass to
a commutative coefficient ring.  We keep track of the filtration
maps, since their compatibility identifies the capacities.

\begin{proposition}
    \label{prop:integral-matijevic-comparison}
    Let $C\subset\bbR^{2n}$ be a non-degenerate strongly convex domain
    containing the origin in its interior, and let $R$ be a commutative
    unital ring.  There is an isomorphism of filtered graded $R$-modules,
    indexed by $L>0$ outside $\sigma(\partial C)$,
    \[
        D_{L,R}^C\colon
        H_*^{S^1}\bigl(\Lambda_C^{<L};R\bigr)
        \xrightarrow{\ \cong\ }
        \SH_{*+n+1}^{S^1,+,<L}(C;R).
    \]
\end{proposition}

\begin{proof}
We lift Matijevi\'c's chain construction
\cite[Theorem~8.2]{matijevic2024positive} to $\bbZ$.
For an $S^1$-invariant Hamiltonian $H$ satisfying its hypotheses,
coherent orientations give a filtered chain isomorphism
\[
    \Theta_{\bbZ}^{S^1,N}\colon
    \operatorname{CM}_{*}^{S^1,N}
    (\psi_{H^*}^{l},\widetilde f_N;\bbZ)
    \xrightarrow{\ \cong\ }
    \operatorname{CF}_{*+n}^{S^1,N}
    (H,\widetilde f_N;\bbZ).
\]
Here $\psi_{H^*}^{l}$ denotes the reduced dual functional, and the map
counts hybrid half-cylinders.
For the non-equivariant hybrid problem, we use the orientation-line
construction of Abbondandolo, Edtmair, and Kang
\cite[Appendix~A.4 and Appendix~B, (B.4)]{abbondandolo2024closed}.
The additional point at a finite Borel stage is that the parameter
equation is independent of the loop.  Thus the linearized half-cylinder
operator, including boundary evaluation, is block upper triangular
\cite[Proposition~3.2 and (4.1)--(4.2)]{matijevic2024positive}.
Scaling its off-diagonal block to zero gives the orientation-line
identification
\[
 o(D)\cong o(D_{\mathrm{Floer}})\otimes o(D_{\mathrm{par}}),
\]
where the two diagonal operators include, respectively, the positive
Fourier trace and the parameter endpoint evaluation.
The Morse boundary condition is imposed by taking the inverse image
of the graph of the unstable manifold.  The determinant exact sequence
for this inverse image tensors the displayed line with the orientation
line of that unstable manifold, exactly as in
\cite[(A.6) and (B.4)]{abbondandolo2024closed}.

Use $S^1$-transported trivializations along the asymptotic circles and
the standard orientation of the circle to orient the free quotient,
as in \cite[Remark~4.5]{bourgeois2009gysin}.
The block deformation and boundary evaluation are equivariant, so
the preceding identification descends to this quotient.
It is compatible with gluing: the determinant exact sequences are
natural under gluing, and at the split operator the signs are those
of the non-equivariant hybrid problem and ordinary Morse theory.
The signed boundary calculation of
\cite[Appendix~B]{abbondandolo2024closed} therefore gives the chain-map
identity for $\Theta_{\bbZ}^{S^1,N}$.

The action-energy inequality implies that $\Theta_{\bbZ}^{S^1,N}$ is upper
triangular with respect to the action filtration.  Since its diagonal entries
are the signed counts of the constant solutions, they equal $\pm1$.
Consequently, $\Theta_{\bbZ}^{S^1,N}$ is a filtered chain isomorphism over
$\bbZ$, and tensoring this chain isomorphism and its inverse with $R$ yields
a filtered chain isomorphism over any commutative ring $R$. Moreover, applying the same orientation
conventions to the continuation and chain-homotopy moduli spaces shows that
these isomorphisms commute up to chain homotopy with the continuation maps.
This also covers the maps between Borel stages in
\cite[Theorem~8.2(2)]{matijevic2024positive}: under the standard inclusions,
the added parameter directions have their complex orientation, and the
same determinant-line identifications restrict to those at the smaller stage.

With these signed maps and homotopies in place, the remainder of
Matijevi\'c's proof \cite[Theorem~C.2, \S11]{matijevic2024positive}
applies unchanged.  For Hamiltonians $H=\varphi\circ H_C$ adapted to the
level $L$, it compares the Clarke sublevel with relative sublevels of the
Hamiltonian dual functional, then uses the hybrid comparison to identify
their homology with positive Floer homology.  Passing to the Borel and
continuation limits gives $D_{L,R}^C$, compatible with the filtration maps.
The first comparison uses equivariant deformations and the long exact
sequence of a pair and already holds with arbitrary coefficients.
In the finite-dimensional reduction, transport the
orientations of the unstable manifolds through the graph identifications
$\Gamma_{l_2,l_1}$ of \cite[Lemma~10.7, (10.9)]{matijevic2024positive};
the additional convex directions introduce no unstable directions, so
these identifications also give the independence of $l$ over $\bbZ$.
\end{proof}

\begin{corollary}
\label{cor:Clarke-Floer-capacities}
For every compact convex body $C\subset\bbR^{2n}$ containing the
origin in its interior, every commutative unital ring $R$, and every
$k\geq1$,
\begin{equation}
    c_k^{R,\Psi}(C)=c_k^R(C).
    \label{eq:Clarke-Floer-capacities}
\end{equation}
\end{corollary}

\begin{proof}
First suppose that $C$ is non-degenerate and strongly convex.
The compatible maps $D_{L,R}^C$ induce an isomorphism $D_{\infty,R}^C$
on direct limits, as in
\cite[Theorem~C.2 and Remark~11.1]{matijevic2024positive}.
The Borel projection factors as
\[
 (\Lambda_C^{<L})_{hS^1}\longrightarrow\Lambda_{hS^1}
 \longrightarrow BS^1.
\]
Since $\Lambda$ is contractible, the second map is a homotopy
equivalence.  Thus $(p_L^C)_*$ is the Clarke filtered-to-limit map
with its target identified with $H_*(BS^1;R)$.

On the Floer side, the unfiltered connecting map identifies positive
homology with shifted low-action homology, since full symplectic
homology vanishes.  Naturality of the low-action/positive exact
sequence identifies the filtered connecting map with the
filtered-to-limit map followed by this identification.  In the
normalization of \cref{eq:coefficient-capacity-definition}, we therefore
obtain the commutative square
\[
 \xymatrix{
 H_*^{S^1}(\Lambda_C^{<L};R)
   \ar[r]^{D_{L,R}^C}\ar[d]_{(p_L^C)_*}
 &
 \SH_{*+n+1}^{S^1,+,<L}(C;R)\ar[d]^{\rho_{L,R}^C}
 \\
 H_*(BS^1;R)\ar[r]_{d_R^C}^{\cong}
 & H_*(BS^1;R).
 }
\]
Here $d_R^C$ is induced by $D_{\infty,R}^C$ under these identifications.
Over $\bbZ$, it is multiplication by $\pm1$ in each relevant degree.
After base change, an $R$-submodule contains $\beta_{k-1}$ if and only
if it contains this unit multiple.  The two image conditions therefore
agree at every non-spectral level.  Such levels suffice for the
capacity infima, proving the equality for $C$.

To extend the equality to every convex body by convex approximation,
it suffices to check inclusion monotonicity and conformality for
$c_k^{R,\Psi}$.  The support-function relations $s_C\leq s_D$ for
$C\subset D$ and $s_{rC}=rs_C$, together with
\eqref{eq:Clarke-support-dual}, give
$\Psi_C\leq\Psi_D$ and $\Psi_{rC}=r^2\Psi_C$.
The induced inclusions of Clarke sublevels commute with projection to
$BS^1$, proving the required properties.
\end{proof}

For ellipsoids, the comparison turns the homological image calculated in
\cref{prop:ellipsoid-Borel-sphere} into a capacity formula.

\begin{proof}[Proof of \cref{thm:cyclic-ellipsoid-capacities}]
Take $L$ outside the ellipsoid action spectrum. By
\cref{prop:ellipsoid-Borel-sphere}, the image of the strict sublevel
at $L$ is the kernel of cap product with
\[
 e_{\bbZ/p^e}(\xi_L)=\Delta_a(L)\,u^{N_a(L)},
\]
as in \eqref{eq:ellipsoid-Euler-class}.
Its cap product with $\beta_{k-1}$ vanishes precisely when
$N_a(L)\geq k$ or $p^e\mid\Delta_a(L)$. Taking the infimum over such
levels therefore gives
\[
 c_k^{\bbZ/p^e,\Psi}(E(a))
 =\min\{c_k^{\bbQ}(E(a)),\tau_{p^e}(a)\},
\]
by \eqref{eq:rational-ellipsoid-count} and the definition of $\tau_{p^e}$.
The comparison
\eqref{eq:Clarke-Floer-capacities} gives the Floer-capacity formula.
\end{proof}

\section{Torsion obstructions from two glued spheres}
\label{sec:243-domains}

We prove \cref{thm:243-obstruction} for the source
\[
 \Omega_X=\left\{(x,y)\in\bbR_{\geq0}^2\ \middle|\
              x\leq1,\quad 6x+5y\leq7\right\},
 \qquad X=X_{\Omega_X},\qquad Y=E(1,3).
\]
The proof computes the capacity values
\begin{equation}
 \label{eq:243-capacity-values}
 c_{15}^{\bbZ/243}(X)=\frac{48}{5}>9=c_{15}^{\bbZ/243}(Y),
\end{equation}
while verifying all the rational and ECH inequalities.  We first explain
how the two requirements guide the choice of source.  The calculation
then proceeds through its Clarke sublevel, and at the end of the
section we discuss the patterns it suggests for other coefficients.

\begin{figure}[tbp]
 \centering
 % Positive polar: forced region and chosen source.
\begin{tikzpicture}[x=4.3cm,y=4.3cm,>=stealth,
  every node/.style={font=\small}]
 \fill[torsionLevelEight] (0,0)--(1,0)--(2/3,2/3)--(0,2/3)--cycle;
 \draw[figureGuideColor,densely dotted]
      (0,4/3)--(4/3,0);
 \draw[sourceBoundary,thick]
      (0,0)--(1,0)--(6/7,5/7)--(0,5/7)--cycle;
 \draw[figureGuideColor,dashed]
      (1,0)--(2/3,2/3)--(0,2/3);
 \draw[->,figureAxisColor] (0,0)--(1.52,0) node[right] {$r$};
 \draw[->,figureAxisColor] (0,0)--(0,1.5) node[above] {$s$};
 \node[below left] at (0,0) {$0$};
 \node at (.30,.30) {$P$};
 \node[sourceBoundary] at (.28,.87) {$\Omega_X^*$};
 \node[anchor=west] at (.12,1.38) {$r+s=4/3$};
 % Crosses mark the four forbidden points on r+s=4/3.
 \foreach \x/\y in {1.333333/0,1/.333333,.333333/1,0/1.333333} {
   \draw[targetBoundary,thick]
      (\x-.018,\y-.018)--(\x+.018,\y+.018)
      (\x-.018,\y+.018)--(\x+.018,\y-.018);
 }
 \fill[figurePointColor] (1,0) circle[radius=1.4pt];
 \fill[figurePointColor] (2/3,2/3) circle[radius=1.6pt];
 \draw[figureGuideColor,thin] (2/3,2/3)--(.96,.87);
 \node[anchor=west] at (.96,.87) {$(2/3,2/3)$};
 \node[below] at (1,0) {$1$};
 \node[below] at (4/3,0) {$4/3$};
 \node[left] at (0,4/3) {$4/3$};
 \node[left] at (0,2/3) {$2/3$};
\end{tikzpicture}
 \caption{Choosing the source polar.  The blue boundary is
 $\Omega_X^*$, and the shaded polygon is
 $P=\operatorname{conv}\{(0,0),(1,0),(2/3,2/3),(0,2/3)\}$.
 The crosses mark the four other points $(r/3,s/3)$ with
 $r,s\in\bbN$ and $r+s=4$.}
 \label{fig:243-source-choice}
\end{figure}
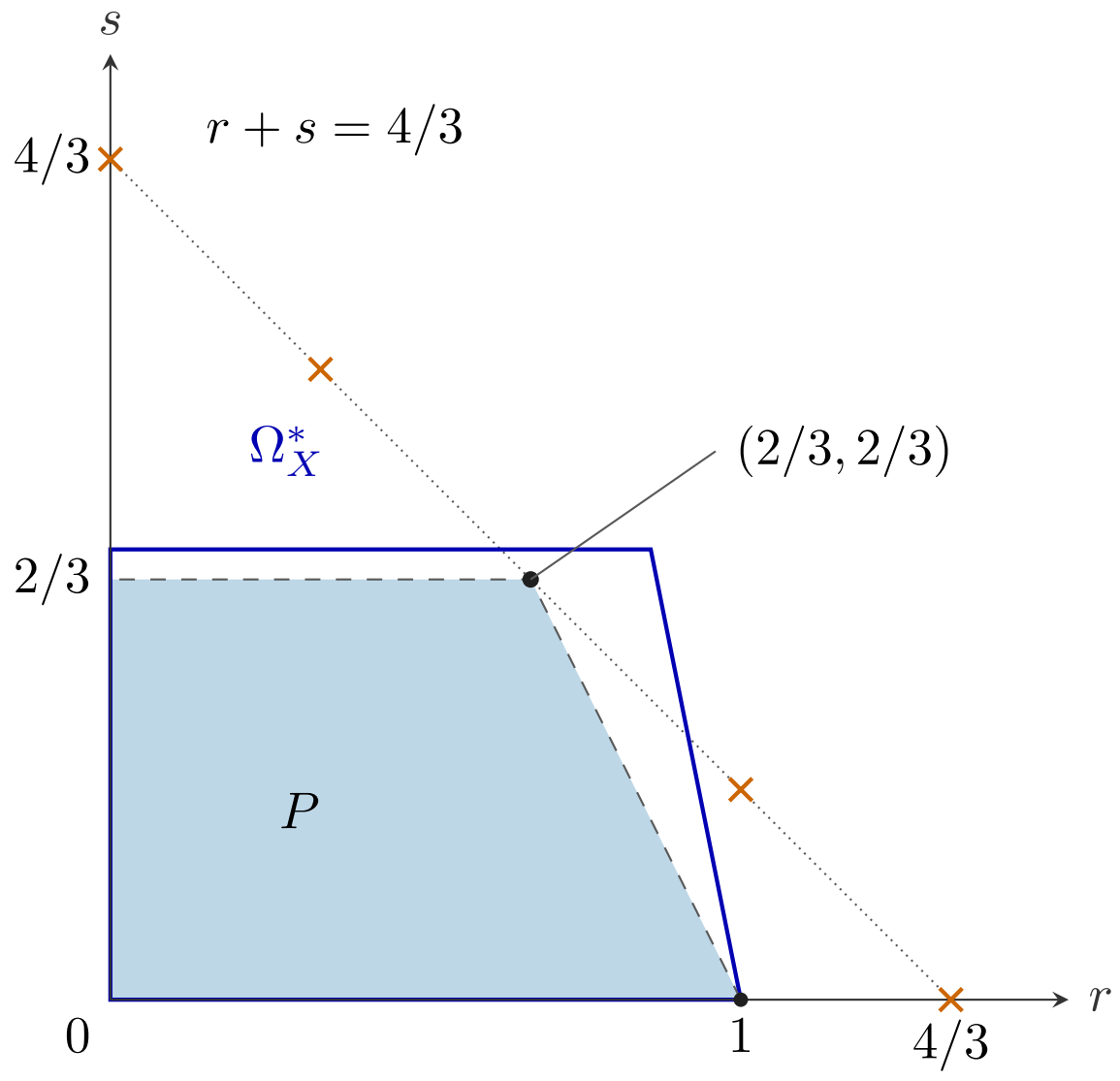

In \cref{fig:243-source-choice}, the two black points provide the
directions needed for the rational comparison, while admitting any
of the crosses would force $T_{243}\leq9$: multiplying its coordinates
by $9$ gives a lattice point with factorial product divisible by
$243$.  The blue boundary shows our choice of polar, containing the
shaded polygon while leaving all four crosses outside.

The support function is
\[
 h_X(r,s)=\max\left\{r+\frac{s}{5},\frac{7s}{5}\right\}.
\]
The two marked directions in \cref{fig:243-source-choice} give
$h_X(1,0)=1$ and $h_X(2,2)=14/5<3$.  The direction $(2,2)$
costs less than the target's $3$ for four units of coordinate sum;
$(1,0)$ supplies the remaining units at cost one each.
Writing $j=4q+t$ with $0\leq t\leq3$, subadditivity gives
\[
 c_j^{\bbQ}(X)
 \leq h_X\bigl(q(2,2)+t(1,0)\bigr)
 \leq \frac{14q}{5}+t
 \leq3q+t=c_j^{\bbQ}(Y).
\]
The ECH inequalities follow from the convex lattice path formula
\cite[Proposition~5.6]{hutchings2016beyond}.
For an integer action level $L=3q+r$, with $0\leq r\leq2$, the path
$(0,2q+r)\to(2q,r)\to(2q,0)$ has $X$-length
$2q\,h_X(1,1)+r\,h_X(1,0)=14q/5+r\leq L$.
The region between this path and the coordinate axes contains
$(2q+1)(q+r+1)$ lattice points, including its boundary, whereas
there are $(q+1)(3q+2r+2)/2$ pairs $(a,b)\in\bbN^2$ with $a+3b\leq L$.
The former count exceeds the latter by $q(q+2r+1)/2\geq0$.
Since the latter counts the ECH capacities of $Y$ at most $L$,
including $c_0^{\mathrm{ECH}}$, the lattice path formula and
monotonicity in the index give
$c_j^{\mathrm{ECH}}(X)\leq c_j^{\mathrm{ECH}}(Y)$ for every $j$.

For this source, the sublevel calculation requires only two maximal
collections of Fourier modes.

For the target torsion capacity, the first level at which
$\lfloor L\rfloor!\lfloor L/3\rfloor!$ is divisible by $243$ is $9$:
\[
 \nu_3(9!\,3!)=5,\qquad \nu_3(8!\,2!)=2.
\]
Since $c_{15}^{\bbQ}(Y)=12$, the ellipsoid formula gives
$c_{15}^{\bbZ/243}(Y)=9$.  For the source, the vector $(9,3)$
gives the upper bound
\[
 c_{15}^{\bbZ/243}(X)\leq h_X(9,3)=\frac{48}{5}.
\]

Put $R=\bbZ/243$.  To prove the lower bound, we will show that
$\beta_{14}$ is absent from the Borel image at the strict level
$48/5$.  The proof has three parts.  A reparametrization replaces
pointwise minimization in the dual Hamiltonian by a single parameter
along each loop.  Fourier modes then identify the sublevel with two
spheres glued along a coordinate subsphere.  Finally, the Thom
cofibers of this diagram determine the required homological image.

\subsection{The analytic step}
\label{subsec:243-analytic-step}

We want to reduce the Clarke sublevel to quadratic sublevels, which
can then be analyzed by Fourier modes.  We first express the dual
Hamiltonian as a minimum of quadratic forms; reparametrization will
then allow us to choose a single minimizing parameter along each loop.

The Hamiltonian of the source is
\[
 H_X(z)=\max\left\{\pi|z_1|^2,
                   \frac{\pi}{7}(6|z_1|^2+5|z_2|^2)\right\}.
\]
Consider the convex combinations of these two quadratic forms,
\[
 H_\lambda(z)=\lambda\pi|z_1|^2
       +(1-\lambda)\frac{\pi}{7}(6|z_1|^2+5|z_2|^2),
 \qquad 0\leq\lambda\leq1.
\]
Their Legendre transforms are the quadratic forms
\[
 q_\lambda(p)=H_\lambda^*(p)=\frac{7}{4\pi}
       \left(\frac{|p_1|^2}{6+\lambda}
                    +\frac{|p_2|^2}{5(1-\lambda)}\right),
 \qquad 0\leq\lambda<1.
\]
At $\lambda=1$, the second summand is zero when $p_2=0$ and
$+\infty$ otherwise.  The associated quadratic functionals on loops are
\[
 Q_\lambda(x)=\int_{\mathbb T}q_\lambda(-J_0\dot x)\,dt.
\]
The next lemma recovers $H_X^*$ from this family.  It is the minimax
identity for the Legendre transform of a maximum of convex functions,
applied here to the two defining quadratic forms.

\begin{lemma}
\label{lem:243-dual}
The dual Hamiltonian satisfies
\[
 H_X^*(p)=\min_{\lambda\in[0,1]}q_\lambda(p).
\]
\end{lemma}

\begin{proof}
Since $H_X=\max_{\lambda\in[0,1]}H_\lambda$, the definition of the
Legendre transform gives
\[
 H_X^*(p)=\sup_{z\in\bbC^2}\min_{\lambda\in[0,1]}
                \bigl(\langle p,z\rangle-H_\lambda(z)\bigr).
\]
For fixed $p$, the expression in parentheses is continuous, concave
in $z$ and affine in $\lambda$.  Sion's minimax theorem
\cite[corollary to Theorem~3.4, p.~174]{sion1958minimax}, with the
compact parameter interval $[0,1]$, therefore gives
\[
 H_X^*(p)=\inf_{\lambda\in[0,1]}\sup_{z\in\bbC^2}
                \bigl(\langle p,z\rangle-H_\lambda(z)\bigr)
          =\inf_{\lambda\in[0,1]}q_\lambda(p).
\]
The infimum is attained because $q_\lambda(p)$ is lower
semicontinuous in $\lambda$, being a supremum of continuous functions,
and $[0,1]$ is compact.
\end{proof}

The minimizing parameter in this formula depends on the two coordinate
speeds of a loop.  Making those speeds constant will allow us to use
one parameter throughout the integral, while preserving the action.
Reparametrization by arc length is a standard tool for comparing
energy sublevels with spaces of constant-speed loops; see
Goresky--Hingston \cite[Proposition~2.2]{goresky2009loop}, who attribute
this comparison to Anosov.

Here we apply Anosov's equivariant reparametrization separately to the
two coordinates.  The point to check is that the resulting homotopy
also preserves the sublevels of $\mathcal H_X$, which couples the
coordinate speeds.  This follows from the explicit interpolation of
the speeds and convexity, as we show below.

\begin{lemma}
\label{lem:243-reparametrization}
There is a continuous $S^1$-equivariant homotopy
\[
 R_s:H^1_0(\mathbb T,\bbC^2)\longrightarrow H^1_0(\mathbb T,\bbC^2),
 \qquad 0\leq s\leq1,
\]
with $R_0=\mathrm{id}$, preserving each coordinate's symplectic action,
such that both coordinate speeds of $R_1x$ are constant and
\[
 \mathcal H_X(R_sx)\leq\mathcal H_X(x),
 \qquad Q_\lambda(R_sx)\leq Q_\lambda(x)\quad(0\leq\lambda<1).
\]
\end{lemma}
\begin{proof}
For a single coordinate, we use the construction in the proof of
Anosov's equivariant reparametrization theorem
\cite[Theorem~4 and \S7(a)]{anosov1981homotopies}, followed by
subtraction of the loop's mean.  We recall the construction in our
setting to keep track of its effect on the speed and symplectic action.

Let $x\in H^1_0(\mathbb T,\bbC)$, $a=|\dot x|$ and
$\ell=\int_{\mathbb T}a$.  If $\ell=0$, set $R_sx=0$.  Otherwise, let
$b$ be the unique periodic mean-zero primitive of $a/\ell-1$, and set
\[
 f(t)=t+b(t).
\]
Then $f$ is a nondecreasing degree-one map of $\bbR$.  Since $x$ is
constant on every interval collapsed by $f$ and satisfies
$|x(v)-x(t)|\leq\ell(f(v)-f(t))$ for $t\leq v$, there is a unique
periodic $\ell$-Lipschitz curve $\gamma$ of length $\ell$ such that
\[
 x=\gamma\circ f,\qquad |\dot\gamma|=\ell\quad\text{almost everywhere}.
\]

Interpolating this normalized arc-length parameter with the original
parameter, put $h_s=(1-s)f+s\,\mathrm{id}$ and define
\[
 R_sx(t)=\gamma(h_s(t))
          -\int_{\mathbb T}\gamma(h_s(v))\,dv.
\]
For $s>0$, the map $h_s$ has a Lipschitz inverse, and the chain rule gives
\begin{equation}
\label{eq:243-speed}
 |\partial_tR_sx(t)|=(1-s)a(t)+s\ell
 \quad\text{almost everywhere};
\end{equation}
at $s=0$ the same identity follows from $R_0x=x$.  Composition with an
orientation-preserving degree-one reparametrization preserves
$\frac12\int_{\mathbb T}\langle J_0x,\dot x\rangle\,dt$, and this also
holds for the nondecreasing map $f$ by change of variables; subtracting
the mean does not affect the integral.  The mean-zero normalization of
$b$ implements Anosov's choice of the initial point and ensures
equivariance: for $x^\theta(t)=x(t-\theta)$ one has
$f^\theta(t)=f(t-\theta)+\theta$ and
$\gamma^\theta(v)=\gamma(v-\theta)$, so $R_s$ commutes with time shift.

To prove continuity, suppose $x_n\to x$ in $H^1_0$ and $s_n\to s$.
If $\ell>0$, then $f_n\to f$ uniformly.  Equicontinuity and the
identities $x_n=\gamma_n\circ f_n$ imply $\gamma_n\to\gamma$
uniformly, hence $R_{s_n}x_n\to R_sx$ uniformly.  The speed identity
\eqref{eq:243-speed} gives
\[
 \|\partial_tR_sx\|_2^2
 =(1-s)^2\|\dot x\|_2^2+(2s-s^2)\ell^2.
\]
Together with uniform convergence, this identity gives weak
convergence of the derivatives in $L^2$ and convergence of their norms,
hence strong $H^1$-convergence.  If $\ell=0$,
the same identity directly gives $R_{s_n}x_n\to0$ in $H^1_0$.

Apply this construction independently to the two coordinates of
$x=(x_1,x_2)$, using the same homotopy parameter $s$.  Each coordinate's
action is preserved, and \eqref{eq:243-speed} interpolates the vector
of coordinate speeds linearly towards its average.  To see that this
also decreases $\mathcal H_X$, write
$H_X^*(p_1,p_2)=F(|p_1|,|p_2|)$.  The function $F$ is convex on
$\bbR_{\geq0}^2$, being the restriction of the convex function $H_X^*$.
With $a_i=|\dot x_i|$ and $\ell_i=\int_{\mathbb T}a_i$, convexity and
Jensen's inequality give
\[
 \int_{\mathbb T}F((1-s)a_1+s\ell_1,(1-s)a_2+s\ell_2)\,dt
 \leq (1-s)\int_{\mathbb T}F(a_1,a_2)\,dt+sF(\ell_1,\ell_2)
 \leq\mathcal H_X(x).
\]
The same argument for each positive quadratic form $q_\lambda$ proves
the remaining inequality.
\end{proof}

Define
\[
 \overline{\mathcal H}_X(x)=\min_{\lambda\in[0,1]}Q_\lambda(x),
 \qquad
 \overline\Lambda_X^{<L}=\{\cA=1,\ \overline{\mathcal H}_X<L\}.
\]
\begin{corollary}
\label{cor:243-constant-parameter}
The inclusion $\overline\Lambda_X^{<L}\hookrightarrow \Lambda_X^{<L}$ is an
$S^1$-homotopy equivalence for every $L>0$.
\end{corollary}

\begin{proof}
Since the integral of the pointwise minimum is at most the minimum of
the integrals, we have $\overline\Lambda_X^{<L}\subset \Lambda_X^{<L}$.
By \cref{lem:243-reparametrization}, the homotopy $R_s$ preserves both
spaces, and at its endpoint $R_1x$ has constant coordinate speeds, so
$\mathcal H_X(R_1x)=\overline{\mathcal H}_X(R_1x)$ by
\cref{lem:243-dual}.  Thus
$R_1\colon\Lambda_X^{<L}\to\overline\Lambda_X^{<L}$ and the restricted
homotopies give an equivariant homotopy inverse.  Finally, if the
minimum in $\overline{\mathcal H}_X(x)<L$ is attained at $\lambda=1$,
then finiteness of $Q_1(x)$ forces $\dot x_2=0$ almost everywhere, so
$Q_\lambda(x)<L$ persists for $\lambda<1$ sufficiently close to $1$.
\end{proof}

Consequently, \cref{cor:243-constant-parameter} replaces the original
sublevel $\Lambda_X^{<L}$ up to $S^1$-homotopy equivalence by the
union over $\lambda\in[0,1)$ of the quadratic sublevels
$\{x\in\Lambda\mid Q_\lambda(x)<L\}$, each of which separates into
Fourier modes.

\subsection{The sphere diagram}

Set $L_0=48/5$.  At a fixed parameter, damping the modes with
nonnegative coefficient in $Q_\lambda-L_0\cA$ leaves a representation
sphere.  For our source, these negative spaces are all contained
in two such representations.

For $r,s\in\bbN$, write
\[
 V_{r,s}=\bigoplus_{j=1}^r\bbC_j\ \oplus\
                 \bigoplus_{j=1}^s\bbC_j.
\]
In the left panel of \cref{fig:prime-polar-rectangles}, the steep
boundary allows the height to jump from $2$ to $6$ when the width
decreases from $9$ to $8$.  These give the two maximal rectangles
$[0,r]\times[0,s]$ with integer corners in $\{h_X<L_0\}$.
Their intersection has corner $(8,2)$, while the boundary point
$(9,3)$ is excluded.  Each corner records the number of positive
Fourier modes retained in the two coordinates.

As in the ellipsoid calculation, we realize its sphere as
$\cA^{-1}(1)\cap V_{r,s}$, using the diagonal rescaling there to
identify it with the unit sphere.  Set
\[
 K=S(V_{9,2})\cup_{S(V_{8,2})}S(V_{8,6}).
\]
The maps from the common subsphere are the standard coordinate
inclusions.  They are equivariant cofibrations, so this union also
realizes their homotopy pushout.

\begin{lemma}
\label{lem:243-cover}
The inclusion induces an $S^1$-equivariant homotopy equivalence
\[
 K\xrightarrow{\simeq}\Lambda_X^{<48/5}.
\]
\end{lemma}

\begin{proof}
Let $F_\lambda=Q_\lambda-L_0\cA$.  On the $j$-th Fourier mode its
coefficients in the two coordinates are
\[
 \pi j\left(\frac{7j}{6+\lambda}-\frac{48}{5}\right),
 \qquad
 \pi j\left(\frac{7j}{5(1-\lambda)}-\frac{48}{5}\right).
\]
Negative Fourier modes have positive coefficients.  Positive modes
are negative precisely when
\[
 j<\frac{48(6+\lambda)}{35}
 \quad\text{or}\quad
 j<\frac{48(1-\lambda)}7,
\]
respectively.  At $\lambda_*=9/16$, these cutoffs are $9$ and $3$.

For an interval $I\subset[0,1)$, let
\[
 U_I=\{x\in\Lambda\mid F_\lambda(x)<0
                         \text{ for some }\lambda\in I\}.
\]
Each $U_I$ is open, being a union of the open sets $\{F_\lambda<0\}$.
The following invariant open sets cover
$\overline\Lambda_X^{<L_0}$:
\[
 \begin{array}{c|c|c}
 I&\text{negative modes lie in}&
             \text{all these modes are negative at}\\ \hline
 [0,9/16]&V_{8,6}&0\\
 {[9/16,1)}&V_{9,2}&3/5
 \end{array}
\]
For $x\in U_I$, choose a parameter witnessing $F_\lambda(x)<0$.
Damping the modes outside the indicated representation does not
increase $F_\lambda$.  This inequality ensures that the action
remains positive, so we can normalize it to one throughout the
homotopy.  The damping uses the fixed Fourier projection onto the
representation and is independent of the witnessing parameter.
The last column shows that the entire sphere of that
representation belongs to $U_I$.  Thus each sphere inclusion is
an equivariant homotopy equivalence.

For each loop, $F_\lambda(x)$ is convex in $\lambda$, so the set
of parameters for which it is negative is an interval.  It follows that
\[
 U_{[0,9/16]}\cap U_{[9/16,1)}
       =\{x\in\Lambda\mid F_{9/16}(x)<0\}.
\]
At this parameter, the negative space is $V_{8,2}$.
The same damping identifies the intersection with $S(V_{8,2})$.
The sphere inclusions therefore give an objectwise equivariant
homotopy equivalence between the sphere diagram and the open-cover
diagram. Passing to homotopy pushouts, we obtain an equivariant
homotopy equivalence
$K\hookrightarrow\overline\Lambda_X^{<L_0}$.
Together with \cref{cor:243-constant-parameter}, this proves the claim.
\end{proof}

In particular, the induced equivalence of Borel constructions is
over $BS^1$.  We can therefore compute the capacity using the
projection $K_{hS^1}\to BS^1$.

\subsection{The algebraic step}

The Euler classes of the two spheres record the covering
multiplicities.  To account for the gluing, we use the Thom cofiber
sequences simultaneously, keeping track of their common subsphere.

\begin{lemma}
\label{lem:243-image}
With $R=\bbZ/243$, the Borel projection has image
\[
 \im\bigl(H_{28}(K_{hS^1};R)\longrightarrow H_{28}(BS^1;R)\bigr)
       =3R\beta_{14}.
\]
\end{lemma}

\begin{proof}
For a representation $V$, write
$T(V)=\operatorname{Th}(ES^1\times_{S^1}V)$ and let
$z_V:BS^1_+\to T(V)$ be its zero-section map.  The stable cofiber
sequences
\[
 \Sigma^\infty S(V)_{hS^1,+}\longrightarrow
 \Sigma^\infty BS^1_+\xrightarrow{z_V}\Sigma^\infty T(V)
\]
commute with the coordinate inclusions in the sphere diagram.
Taking their homotopy pushout gives a cofiber sequence
\[
 \Sigma^\infty K_{hS^1,+}\longrightarrow
 \Sigma^\infty BS^1_+\xrightarrow{q}\mathcal T,
\]
where $\mathcal T$ is the homotopy pushout of the Thom spectra.
The middle term remains $\Sigma^\infty BS^1_+$ because the indexing
category has contractible nerve.  Exactness gives
\[
 \im\bigl(H_{28}(K_{hS^1};R)\to H_{28}(BS^1;R)\bigr)=\ker q_*.
\]

The Thom isomorphism for a complex representation of dimension $d$ is
\[
 \widetilde H_j(T(V);R)\cong H_{j-2d}(BS^1;R).
\]
The representations $V_{9,2}$, $V_{8,6}$ and $V_{8,2}$ have dimensions
$11$, $14$ and $10$, respectively.  Their degree-$28$ Thom groups
are therefore $R\beta_3$, $R\beta_0$ and $R\beta_4$.

For a coordinate inclusion $V_{r,s}\subset V_{a,b}$, naturality of
the Thom class identifies the induced map on Thom homology with
\[
 H_{j-2(r+s)}(BS^1;R)
 \xrightarrow{\ e(V_{a,b}/V_{r,s})\frown(-)\ }
 H_{j-2(a+b)}(BS^1;R).
\]
Indeed, the
Thom class of a direct sum is the product of the Thom classes, and
restriction to the first summand replaces the second by its Euler
class \cite[Section~3.2, proof of Proposition~3.13(b)]{hatcherVectorBundles}.
Here that class is
\[
 e(V_{a,b}/V_{r,s})
 =\frac{a!\,b!}{r!\,s!}\,u^{a+b-r-s}.
\]
Since $u^m\frown\beta_j=\beta_{j-m}$, the two maps out of the
overlap multiply by $9$ and $3\cdot4\cdot5\cdot6=360$.
All degree-$27$ Thom groups vanish, so Mayer--Vietoris gives
\[
 R\xrightarrow{M}R^2\longrightarrow H_{28}(\mathcal T;R)
          \longrightarrow0,\qquad
 M=\begin{pmatrix}9\\-360\end{pmatrix}.
\]
The rows correspond to $V_{9,2}$ and $V_{8,6}$; the sign is the
usual $(i_*,-j_*)$ convention for the overlap.

The map $q$ factors through the zero section at either vertex.
For $V_{9,2}$, the induced maps fit into the diagram
\[
 \begin{tikzcd}[column sep=large,row sep=large]
 H_{28}(BS^1;R)
   \arrow[r,"(z_{V_{9,2}})_*"] \arrow[dr,"q_*"']
 &
 \widetilde H_{28}(T(V_{9,2});R) \arrow[d]
 \\
 & H_{28}(\mathcal T;R)\cong\operatorname{coker}M.
 \end{tikzcd}
\]
Under the Thom isomorphism, the horizontal map is cap product with
$9!\,2!\,u^{11}$.  Consequently
\[
 q_*(\beta_{14})=[v]\in\operatorname{coker}M,\qquad
 v=(9!\,2!,0)=(162,0)\quad\text{in }R^2.
\]
The row vector $w=(40,1)$ annihilates $M$ and satisfies
$wv=162\neq0$ in $R$.  Thus $t[v]=0$ implies $162t=0$, or $3\mid t$.
Conversely, $3v=0$.  Hence $\ker q_*=3R\beta_{14}$, as claimed.
\end{proof}

\begin{proof}[Proof of \cref{thm:243-obstruction}]
By \cref{lem:243-cover,lem:243-image}, the class $\beta_{14}$ is
not in the Borel image of $\Lambda_X^{<48/5}$.  Every smaller sublevel
maps through this one, so the Clarke characterization
\eqref{eq:Clarke-Floer-capacities} gives
$c_{15}^{\bbZ/243}(X)\geq48/5$.  The supporting-ellipsoid bound
computed at the beginning of the section gives the reverse inequality.
Together with the target calculation, this proves
\eqref{eq:243-capacity-values}.  Symplectic monotonicity excludes the
embedding, while the rational inequalities were established above.
\end{proof}

\subsection*{The pattern for odd primes}

The calculation for $p=3$ suggests the following pattern for odd primes.
Set
\[
 N_p=\frac{(p-1)(p+3)}2,\qquad X_p=X_{\Omega_p},
\]
where
\[
 \Omega_p=\left\{(x,y)\in\bbR_{\geq0}^2\ \middle|\
       x\leq1,\quad N_px+(p+2)y\leq N_p+1\right\}.
\]
Thus $X_3=X$.  With $k=p^2+2p$, the expected capacity values are
\begin{equation}
 \label{eq:odd-prime-capacity-values}
 c_k^{\bbZ/p^{p+2}}(X_p)
       =p^2+\frac{p}{p+2}>p^2
       =c_k^{\bbZ/p^{p+2}}(E(1,p)),
\end{equation}
while all rational capacity inequalities hold.

The rational comparison follows the same geometric choice as for
$p=3$.  The support function
\[
 h_{\Omega_p}(r,s)=
       \max\left\{r+\frac{s}{p+2},\frac{(N_p+1)s}{p+2}\right\}
\]
assigns cost $1$ to $(1,0)$ and cost $p-1/(p+2)<p$ to
$(p-1,2)$, whose coordinate sum is $p+1$.
Writing $j=(p+1)q+t$ with $0\leq t\leq p$, we therefore have
\[
 c_j^{\bbQ}(X_p)
 \leq h_{\Omega_p}\bigl(q(p-1,2)+t(1,0)\bigr)
 \leq pq+t=c_j^{\bbQ}(E(1,p)).
\]
The point $(p^2,p)$, whose factorial product has $p$-adic valuation
$p+2$, gives the source upper bound
$L_p=p^2+p/(p+2)$.  The ellipsoid formula gives the target value $p^2$.

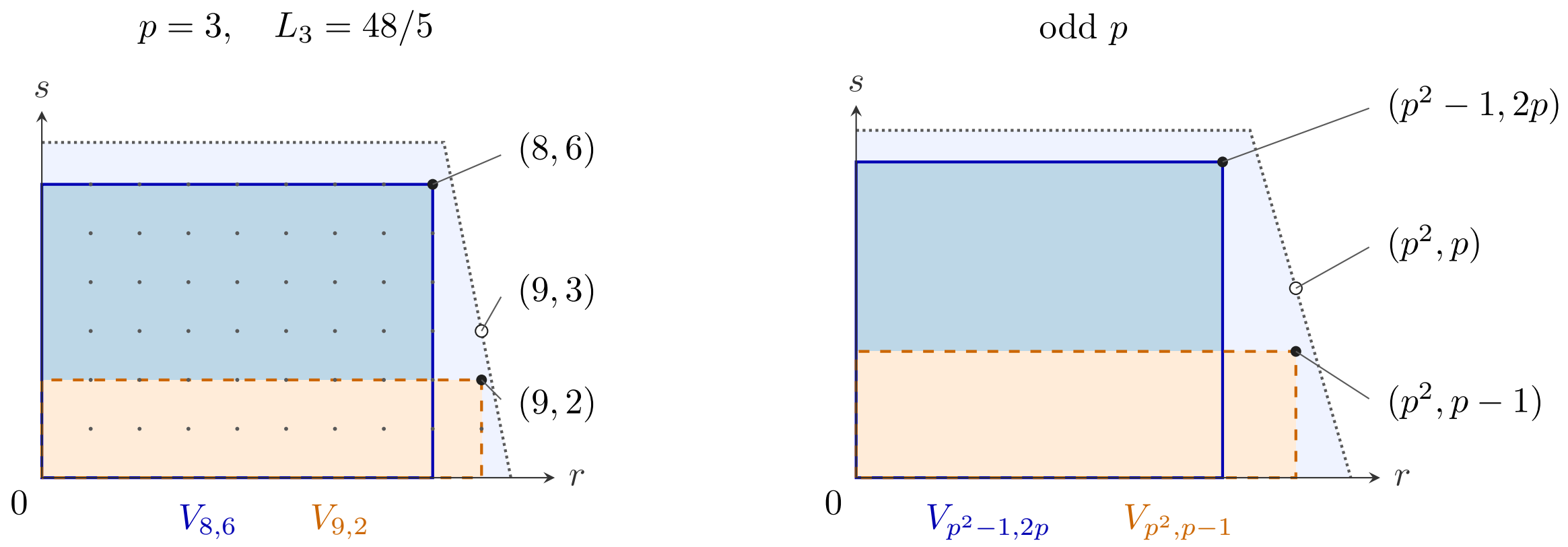
\begin{figure}[H]
 \centering
 % Left: exact p=3 geometry. Right: schematic lattice rectangles for odd p.
\begin{tikzpicture}[>=stealth,every node/.style={font=\small}]
 \begin{scope}[x=.48cm,y=.48cm]
  \node at (5,9.2) {$p=3,\quad L_3=48/5$};
  \fill[torsionLevelSeven] (0,0)--(9.6,0)--(288/35,48/7)--(0,48/7)--cycle;
  \fill[torsionLevelEight] (0,0) rectangle (8,6);
  \fill[targetFill] (0,0) rectangle (9,2);
  \draw[sourceBoundary,thick] (0,0) rectangle (8,6);
  \draw[targetBoundary,thick,dashed] (0,0) rectangle (9,2);
  \draw[figureGuideColor,thick,densely dotted]
       (0,48/7)--(288/35,48/7)--(9.6,0);
  \foreach \r in {1,...,9} {
   \foreach \s in {1,...,6} {
    \ifdim\r pt=9pt
      \ifdim\s pt<3pt
       \fill[figureGuideColor] (\r,\s) circle[radius=.55pt];
      \fi
    \else
      \fill[figureGuideColor] (\r,\s) circle[radius=.55pt];
    \fi
   }
  }
  \draw[->,figureAxisColor] (0,0)--(10.5,0) node[right] {$r$};
  \draw[->,figureAxisColor] (0,0)--(0,7.5) node[above] {$s$};
  \fill[figurePointColor] (8,6) circle[radius=1.5pt];
  \fill[figurePointColor] (9,2) circle[radius=1.5pt];
  \draw[figurePointColor,fill=white] (9,3) circle[radius=1.7pt];
  \draw[figureGuideColor,thin] (8,6)--(9.4,6.6);
  \node[anchor=west] at (9.45,6.7) {$(8,6)$};
  \node[anchor=west] at (9.45,1.5) {$(9,2)$};
  \draw[figureGuideColor,thin] (9,2)--(9.4,1.6);
  \node[anchor=west] at (9.45,3.8) {$(9,3)$};
  \draw[figureGuideColor,thin] (9,3)--(9.4,3.7);

  \node[below left] at (0,0) {$0$};
  \node[sourceBoundary] at (3.4,-.9) {$V_{8,6}$};
  \node[targetBoundary] at (6.1,-.9) {$V_{9,2}$};
 \end{scope}
 \begin{scope}[shift={(8cm,0)},x=.72cm,y=.62cm]
  \node at (3.1,7.1) {odd $p$};
  \fill[torsionLevelSeven] (0,0)--(6.75,0)--(5.375,5.5)--(0,5.5)--cycle;
  \fill[torsionLevelEight] (0,0) rectangle (5,5);
  \fill[targetFill] (0,0) rectangle (6,2);
  \draw[sourceBoundary,thick] (0,0) rectangle (5,5);
  \draw[targetBoundary,thick,dashed] (0,0) rectangle (6,2);
  \draw[figureGuideColor,thick,densely dotted]
       (0,5.5)--(5.375,5.5)--(6.75,0);
  \draw[->,figureAxisColor] (0,0)--(7.15,0) node[right] {$r$};
  \draw[->,figureAxisColor] (0,0)--(0,5.9) node[above] {$s$};
  \fill[figurePointColor] (5,5) circle[radius=1.5pt];
  \fill[figurePointColor] (6,2) circle[radius=1.5pt];
  \draw[figurePointColor,fill=white] (6,3) circle[radius=1.7pt];
  \draw[figureGuideColor,thin] (5,5)--(7,5.85);
  \node[anchor=west] at (7.05,5.9) {$(p^2-1,2p)$};
  \draw[figureGuideColor,thin] (6,2)--(7,1.25);
  \node[anchor=west] at (7.05,1.2) {$(p^2,p-1)$};
  \draw[figureGuideColor,thin] (6,3)--(7,3.65);
  \node[anchor=west] at (7.05,3.7) {$(p^2,p)$};

  \node[below left] at (0,0) {$0$};
  \node[sourceBoundary] at (1.8,-.7) {$V_{p^2-1,2p}$};
  \node[targetBoundary] at (4.4,-.7) {$V_{p^2,p-1}$};
 \end{scope}
\end{tikzpicture}
 \caption{Maximal lattice rectangles in the strict positive-polar
 sublevel: $p=3$ on the left and the suggested odd-prime pattern
 on the right.  Solid blue and dashed orange boundaries distinguish
 the two rectangles, labelled by their Fourier representations.
 The dotted outer boundary is excluded, as is the hollow point on it.
 The right panel is schematic.}
 \label{fig:prime-polar-rectangles}
\end{figure}

The polar picture explains the expected lower bound.
A lattice rectangle with upper-right corner $(r,s)$ records the
representation $V_{r,s}$: it contains the first $r$ positive modes
in one coordinate and the first $s$ in the other.
In the strict polar sublevel
\[
 \{h_{\Omega_p}<L_p\}
 =\left\{(r,s)\in\bbR_{\geq0}^2\ \middle|\
       r+\frac{s}{p+2}<L_p,\quad
       \frac{(N_p+1)s}{p+2}<L_p\right\},
\]
the maximal lattice rectangles have corners $(p^2,p-1)$ and
$(p^2-1,2p)$; their intersection has corner $(p^2-1,p-1)$.
The point $(p^2,p)$ lies on the excluded boundary.
Thus the same picture as in \cref{fig:prime-polar-rectangles}
suggests the sphere model
\[
 K_p=S(V_{p^2,p-1})
       \cup_{S(V_{p^2-1,p-1})}S(V_{p^2-1,2p})
\]
for $\Lambda_{X_p}^{<L_p}$.  For $p=3$, this is precisely the
model established in \cref{lem:243-cover}.

The complementary Euler coefficients are $p^2$ and
$(2p)!/(p-1)!$, both of valuation $2$, while the zero-section
coefficient $(p^2)!(p-1)!$ has valuation $p+1$.
With $R_p=\bbZ/p^{p+2}$, the Thom calculation for this diagram
gives image $pR_p\beta_{k-1}$ in degree $2k-2$, as in
\cref{lem:243-image}.  The model would therefore exclude
$\beta_{k-1}$ and give the lower bound predicted in
\eqref{eq:odd-prime-capacity-values}.

The expected phenomenon at exponent $e=p+2$ is not restricted to
odd primes, although the two-sphere recipe above is.  For $p=2$,
in two complex dimensions one has
$T_{16}(\Omega)=2c_3^{\bbQ}(X_\Omega)$, so the rational inequalities
against $E(1,2)$ prevent the required separation of thresholds.
An analogous construction with $\bbZ/16$ coefficients uses a
six-dimensional domain and a more elaborate gluing of representation
spheres.  Thus the restriction to odd primes concerns the particular
domains and two-sphere model presented here, rather than the
expected coefficient phenomenon.

\section{Stabilization}
\label{sec:stabilization-sketch}

We prove that taking products with disks of sufficiently large area
preserves the capacities of convex bodies with arbitrary coefficients,
as stated in \cref{thm:capacity-disk-stabilization}.
We then apply this result to the torsion obstruction of
\cref{sec:243-domains} to obtain the sharp stabilized embedding threshold.

For rational capacities of convex toric domains, invariance under
adding a disk of sufficiently large area follows from the product
formula of Gutt--Hutchings
\cite[Theorem~1.6 and Remark~1.10]{gutt2018symplectic}.
A related principle appears in Siegel's stabilization theorem for
his SFT capacities: the data defining the capacities are preserved
below the area of the added disk
\cite[\S6.2.2]{siegel2025higher}.

The variational argument uses the support function of a product
and the splitting of symplectic action, as in the first proof of
Haim-Kislev--Ostrover's product formula for the EHZ capacity
\cite[Proposition~1.5 and \S2]{haimkislev2023remarks}.
To control the whole sublevel, we damp the disk component of each
loop and normalize its action, as in
\cref{thm:ellipsoid-weighted-sphere}.  The product formula below gives
the estimate needed for this deformation.

We use the sublevels and Borel projections of
\cref{sec:variational-comparison}, together with convex approximation
and the comparison \eqref{eq:Clarke-Floer-capacities}.

The stabilization theorem follows from a stronger statement:
below the area of the added disk, the inclusion preserves the
$S^1$-equivariant homotopy type of the Clarke sublevels, compatibly
with the sublevel inclusions.

\begin{lemma}
\label{lem:Clarke-disk-stabilization}
Let $C\subset\bbC^n$ be a compact convex body containing the origin in
its interior, and let $A>0$.  The inclusions $x\mapsto(x,0)$ admit
an $S^1$-equivariant strong deformation retraction of filtered spaces
\[
 \Lambda_{C\times D(A)}^{<L}\searrow\Lambda_C^{<L},
 \qquad 0<L<A.
\]
\end{lemma}

\begin{proof}
Maximizing separately over the two factors gives
\begin{align*}
 s_{C\times D(A)}(p,q)
 &=\sup_{\substack{z\in C\\w\in D(A)}}
       \bigl(\langle p,z\rangle+\langle q,w\rangle\bigr)\\
 &=s_C(p)+\sqrt{A/\pi}\,|q|.
\end{align*}
For $q\neq0$, the disk contribution is attained at the boundary point
$\sqrt{A/\pi}\,q/|q|$ pointing in the direction of $q$, and it vanishes
when $q=0$.  By \eqref{eq:Clarke-support-dual}, we obtain
\begin{equation}
 \label{eq:Clarke-product-dual}
 \begin{aligned}
 H_{C\times D(A)}^*(p,q)
 &=\frac{s_{C\times D(A)}(p,q)^2}{4}\\
 &=\frac14\left(s_C(p)+\sqrt{A/\pi}\,|q|\right)^2\\
 &=\left(\sqrt{H_C^*(p)}
             +\sqrt{\frac{A}{4\pi}}\,|q|\right)^2.
 \end{aligned}
\end{equation}

For a mean-zero loop $(x,y)\in H^1_0(\mathbb T,\bbC^n\times\bbC)$,
put
\[
 E(y)=\frac{1}{4\pi}\int_{\mathbb T}|\dot y|^2\,dt,
 \qquad
 B(x,y)=\sqrt{\frac{A}{4\pi}}
       \int_{\mathbb T}\sqrt{H_C^*(-J_0\dot x)}\,|\dot y|\,dt.
\]
The Fourier formulas of \cref{sec:ellipsoid-Clarke-homotopy} give
\[
 E(y)-\cA(y)
 =\pi\sum_{m\neq0}(m^2-m)|y_m|^2\geq0.
\]
Fix $0<L<A$.  Damping the disk component must be followed by
normalization of the action, which can increase the dual functional.
We therefore work with
\[
 F_L=\mathcal H_{C\times D(A)}-L\cA,
\]
whose negativity is preserved by positive rescaling.  On the
action-one hypersurface, $F_L<0$ is precisely the sublevel condition.
By \eqref{eq:Clarke-product-dual}, for $0\leq t\leq1$,
\begin{equation}
 \label{eq:Clarke-product-damping}
 \begin{split}
 F_L(x,ty)
 &=\mathcal H_C(x)-L\cA(x)+2tB(x,y)\\
 &\hspace{2em}+t^2\bigl(AE(y)-L\cA(y)\bigr).
 \end{split}
\end{equation}
Since $B(x,y)\geq0$ and
\[
 AE(y)-L\cA(y)\geq(A-L)E(y)\geq0,
\]
damping the disk component does not increase $F_L$:
\[
 F_L(x,ty)\leq F_L(x,y)<0
 \qquad\text{if }(x,y)\in\Lambda_{C\times D(A)}^{<L}.
\]
The same inequality also ensures that the action remains positive:
\[
 L\cA(x,ty)>\mathcal H_{C\times D(A)}(x,ty)\geq0.
\]
We can therefore normalize throughout the deformation by setting
\[
 R_s(x,y)=
 \frac{(x,(1-s)y)}
      {\sqrt{\cA(x)+(1-s)^2\cA(y)}},
 \qquad 0\leq s\leq1;
\]
positive two-homogeneity and \eqref{eq:Clarke-product-damping} show
that $R_s$ takes values in $\Lambda_{C\times D(A)}^{<L}$.

Because coordinate scaling and the action functional are continuous on
$H^1_0$ and the denominator is positive, $R_s$ is
jointly continuous in the $H^1$ topology.  It commutes with time shift,
starts at the identity, fixes every loop $(x,0)$ of action one, and at
$s=1$ maps $\Lambda_{C\times D(A)}^{<L}$ into $\Lambda_C^{<L}$ since
\[
 \mathcal H_{C\times D(A)}(x,0)=\mathcal H_C(x).
\]
Finally, the formula
for $R_s$ is independent of $L$, which proves compatibility with the
sublevel inclusions.
\end{proof}

\begin{proof}[Proof of \cref{thm:capacity-disk-stabilization}]
For one disk factor, the inclusion $i_L:x\mapsto(x,0)$ induces a
homotopy equivalence on Borel constructions for $L<A$, and
\[
 p_L^{C\times D(A)}\circ(i_L)_{hS^1}=p_L^C.
\]
The two Borel projections therefore have the same homological image
in $H_*(BS^1;R)$ at every such level.  Taking the infima in
\eqref{eq:Clarke-capacity-definition} gives
\[
 c_k^{R,\Psi}(C\times D(A))=c_k^{R,\Psi}(C)
 \qquad\text{whenever }c_k^{R,\Psi}(C)<A.
\]
The comparison
\eqref{eq:Clarke-Floer-capacities}, including its extension to convex
bodies, gives the assertion for one factor.  Iteration proves the result.
\end{proof}

We now apply the stabilization theorem to the capacity values
computed for $X$ and $E(1,3)$.

\begin{proof}[Proof of \cref{cor:stabilized-243}]
Put $R=\bbZ/243$.  The calculation
\eqref{eq:243-capacity-values} gives
\[
 c_{15}^R(X)=48/5,\qquad c_{15}^R(E(1,3))=9.
\]

Suppose that $\varphi$ is an embedding as in the statement.
Fix $0<r<1$ and $A>(48/5)r$, and set
\[
 K=(\sqrt r\,X)\times D(A)^N.
\]
This is a compact subset of the source of $\varphi$.
Choose $B>9\lambda$ large enough that
\[
 \varphi(K)\subset\operatorname{int}(Q),
 \qquad Q=E(\lambda,3\lambda)\times D(B)^N.
\]
By monotonicity and \cref{thm:capacity-disk-stabilization},
\[
 \frac{48}{5}r=c_{15}^R(K)\leq c_{15}^R(Q)=9\lambda.
\]
Thus $\lambda\geq16/15$ by letting $r\uparrow1$.

Conversely, the inequalities defining $\Omega_X$ imply
\[
 3x+y=\frac95x+\frac15(6x+5y)
       \leq\frac95+\frac75=\frac{16}{5}.
\]
Consequently $X\subset E(16/15,16/5)$, with strict inequality for
points of $\operatorname{int}(X)$.  Inclusion, followed by the
identity on $\bbC^N$, gives the required embedding for every
$\lambda\geq16/15$.
\end{proof}

For the family discussed at the end of \cref{sec:243-domains},
the expected capacity values would give the sharp stabilized threshold
\[
 \lambda_p=1+\frac1{p(p+2)}.
\]
Indeed, the support function gives
$h_{\Omega_p}(p,1)=p+1/(p+2)$, so the inclusion
$X_p\subset E(\lambda_p,p\lambda_p)$ realizes this dilation.
The expected source capacity would give the matching lower bound
after stabilization by the argument above.

\printbibliography

\end{document}